\documentclass[a4paper,fleqn]{cas-sc}
\usepackage{subcaption}
\usepackage[authoryear]{natbib}
\def\tsc#1{\csdef{#1}{\textsc{\lowercase{#1}}\xspace}}
\tsc{WGM}
\tsc{QE}
\usepackage{physics}        % abs and norm

\usepackage{algorithm}
\usepackage{algorithmicx}
\usepackage{algpseudocode}

\newcommand{\dispx}{u_{x_{1:t}}}
\newcommand{\dispy}{u_{y_{1:t}}}
\newcommand{\avgstress}{\bar{\sigma}_{1:t}}
\newcommand{\avgstrain}{\bar{\varepsilon}_{1:t}}

\begin{document}
\let\WriteBookmarks\relax
\def\floatpagepagefraction{1}
\def\textpagefraction{.001}

\shorttitle{Cooperative data-driven modeling}    

\shortauthors{A Dekhovich et~al.}  

\title[mode=title]{Cooperative data-driven modeling}


%\tnotetext[1]{Title footnote text} 


\author[1]{Aleksandr Dekhovich}


%\fnmark[<footnote mark no>]

%\ead{<email address>}

%\ead[url]{<URL>}

%\credit{<Credit authorship details>}

\affiliation[1]{organization={Department of Material Science and Engineering, Delft University of Technology},
            addressline={Mekelweg 2}, 
            city={Delft},
%          citysep={}, % Uncomment if no comma needed between city and postcode
            postcode={2628 CD}, 
            %state={},
            country={The Netherlands}}

\author[2]{O. Taylan Turan}
% Footnote of the second author

% Email id of the second author
%\ead{}



\affiliation[2]{organization={Pattern Recognition and Bioinformatics Laboratory, Delft University of Technology},
            addressline={Van Mourik Broekmanweg 6}, 
            city={Delft},
%          citysep={}, % Uncomment if no comma needed between city and postcode
            postcode={2628 XE}, 
            %state={},
            country={The Netherlands}}

\author[1]{Jiaxiang Yi}

\author[3]{Miguel A. Bessa}
% Footnote of the third author
%\fnmark[3]

%Email id of the third author
\ead{miguel_bessa@brown.edu}

% URL of the third author


\affiliation[3]{organization={School of Engineering, Brown University},
            addressline={184 Hope St.}, 
            city={Providence},
%          citysep={}, % Uncomment if no comma needed between city and postcode
            postcode={RI 02912}, 
            %state={},
            country={USA}}

\cormark[1]
\cortext[1]{Corresponding author}

%\fntext[3]{}


\begin{abstract}
Data-driven modeling in mechanics is evolving rapidly based on recent machine learning advances, especially on artificial neural networks. As the field matures, new data and models created by different groups become available, opening possibilities for cooperative modeling. However, artificial neural networks suffer from catastrophic forgetting, i.e. they forget how to perform an old task when trained on a new one. This hinders cooperation because adapting an existing model for a new task affects the performance on a previous task trained by someone else. The authors developed a continual learning method\footnote{The code implementation and data are available at \url{https://github.com/bessagroup/CDDM}.} that addresses this issue, applying it here for the first time to solid mechanics. In particular, the method is applied to recurrent neural networks to predict history-dependent plasticity behavior, although it can be used on any other architecture (feedforward, convolutional, etc.) and to predict other phenomena. This work intends to spawn future developments on continual learning that will foster cooperative strategies among the mechanics community to solve increasingly challenging problems. We show that the chosen continual learning strategy can sequentially learn several constitutive laws without forgetting them, using less data to achieve the same error as standard (non-cooperative) training of one law per model. 
\end{abstract}


%\item Cooperative data-driven modeling (CDDM) framework based on continual learning ideas is proposed.
%\item Continual learning strategy is used for the first time in the context of constitutive law modeling.
%\item CDDM can sequentially learn new constitutive laws without forgetting all previous ones.
%\item By transferring previously learned knowledge, CDDM is able to learn new constitutive laws using fewer data points compared to standard training.

% Keywords
\begin{keywords}
 \sep data-driven modeling \sep continual learning \sep transfer learning \sep plasticity
\end{keywords}

\maketitle

% Main text
\section{Introduction}\label{sec:intro}

Machine learning permeated almost every scientific discipline \citep{shanmuganathan2016artificial, wuest2016machine, schmidt2019recent}, and Solid Mechanics is no exception \citep{bessa2017, capuano2019smart, thakolkaran2022nn}. With all their merits and flaws \citep{karniadakis2021physics}, these algorithms provide a means to understand large datasets, finding patterns and modeling behavior where analytical solutions are challenging to obtain or not accurate enough. This work introduces the concept of cooperative data-driven modeling by highlighting the importance of continual or lifelong learning and exemplifying it in Solid Mechanics. Without loss of generality, the examples provided in this article pertain to using neural networks to create constitutive models from synthetic data \cite{bessa2017}, but the proposed strategy is based on a general method introduced by the authors in the Computer Science community \cite{dekhovich2022continual}, so it is applicable to many other fields that can also benefit from cooperative data-driven modeling.

For readers unfamiliar with the field of using machine learning to learn constitutive models of materials, we provide a short review of the topic. Using neural networks to describe constitutive material behavior was first proposed decades ago by Ghaboussi et al. \cite{ghaboussi1991a} using simple experimental data. However, advances in numerical modeling and the ability to create large synthetic datasets has led to a new era of data-driven modeling initiated in \cite{bessa2017} that is based on fast analysis of representative volume elements of materials. Since then, there has been rapid progress in the field, first by considering similar architectures \cite{ibanez2018a,jones2018a,nguyen2018a,rocha2021a}, then by considering deep learning strategies to characterize more complex behavior \cite{mozaffar2019,ghavamian2019a,vlassis2020a,abueidda2021a,fuhg2022a,liu2022a}, and recently including physical constraints \cite{masi2021a,asad2022a}. In particular, the first work to propose the use of recurrent neural networks for plasticity modeling \cite{mozaffar2019} showed that these architectures could learn the path- and time-dependency behavior of materials. Soon after, several research groups proposed new neural network architectures and solved increasingly complex plasticity problems \citep{zhang2020using, abueidda2021a, saidi2022deep,bonatti2022cp}. A similar trend is ongoing in other fields within and outside Mechanics \citep{liu2019deep, peng2021multiscale, khalil2017learning,dutting2019optimal}.

Simultaneously, the scientific community is experiencing strong incentives to adhere to open science, with vehement support from funding agencies throughout the World to share data and models according to FAIR principles (Findable, Accessible, Interoperable and Reusable) \citep{wilkinson2016fair,draxl2018nomad,jacobsen2020fair}. There is also a clear need for end users to reuse these models and data. Nevertheless, there is a serious issue that obstructs the synergistic use of machine learning models by the community. Artificial neural networks, unlike biological neural networks, suffer from catastrophic forgetting \citep{mccloskey1989catastrophic, french1999catastrophic, goodfellow2013empirical}. Human beings when learning a new task, e.g. playing tennis, do not forget how to perform past tasks, e.g. swimming. Unfortunately, artificial neural networks fail at this because they are based on updating their parameters (weights and biases) for the task and data being considered, but this changes the previous configuration obtained for a past task (that led to different values of weights and biases). This catastrophic forgetting has important implications in practice, as illustrated by the following scenario.

Imagine that Team A of scientists collects computational or experimental data about the constitutive behavior of Material A, and then trains an artificial neural network to predict the behavior of that material. In the end, Team A publishes the artificial neural network model and corresponding data according to FAIR principles. Later, if Team B aims to create a model that predicts the behavior of Material B then it faces two options: 1) collect data and train a model from scratch for this new material; or 2) use the model developed by Team A in an attempt to get a better model for Material B and use less data during the training process. If the material behavior for B has some commonality with the one for A, there is an advantage in leveraging the work from Team A. However, the state of the art in the literature is to use transfer learning or meta-learning methods \citep{caruana1994learning,vilalta2002perspective,pan2010a,tan2018survey} to adapt Model A and retrain it for Team B's scenario \citep{liu2019c,lejeune2020a}. Unfortunately, in this case, the new model obtained by Team B is no longer valid for Team A's scenario. Although Team B may create a model that is valid for its purposes, this would not be a truly cooperative effort with Team A because a general model valid for both scenarios would not be obtained. This represents a significant challenge to cooperative data-driven modeling because it discourages different groups from working towards a common model, ultimately leading to many independent models. Note that the problem gets worse as more tasks accumulate (more materials and more teams).

A new branch of machine learning called continual or lifelong learning \citep{thrun1998, zenke2017continual, parisi2019continual} is recently opening new avenues to address the catastrophic forgetting issue. Despite being at an early stage, we believe that addressing this limitation will unlock a new era of cooperative data-driven modeling traversing all fields of application. This year, the authors proposed a new continual learning algorithm \citep{dekhovich2022continual} and applied it to various standard computer science datasets to demonstrate the best performance to date in the challenging class-incremental learning scenario when compared with state-of-the-art methods. Further developments are needed, but these algorithms represent an essential step towards democratizing cooperative data-driven modeling. Here, the continual learning method is applied to a new architecture suitable for plasticity modeling, demonstrating its benefit and motivating future research in this nascent field.

\section{Methodology}

Continual or lifelong learning \citep{thrun1998} aims to learn tasks in a sequence while only having access to the data of the current task. The main challenge in this setting is catastrophic forgetting \citep{french1999catastrophic, mccloskey1989catastrophic}, i.e. the phenomenon of forgetting how to solve previous tasks while learning a new one. This problem has become highly significant for deep neural networks in computer vision and natural language processing problems in recent years \citep{zenke2017continual, aljundi2018memory, biesialska2020continual, sokar2021spacenet}. Fundamentally, continual learning tasks are divided into two categories \citep{masana2020class, wortsman2020supermasks}: in the first one, the task to be solved during inference is known, i.e. task-ID is given, while in the second it is unknown (task-ID is not given). The former case is more challenging than the latter\footnote{For readers unfamiliar with the concept of not knowing the task-ID, think about the previous example of one of the tasks being 'swimming' and the other being 'playing tennis'. If the task-ID is not known during inference, then understanding that the task being executed is 'playing tennis' has to happen from context (from the data). Otherwise, if the task-ID is given there is no need to infer from context which task is being considered.}. 

Continual learning algorithms are often classified as regularization, rehearsal or architectural methods \citep{masana2020class}. Regularization-based approaches \citep{zenke2017continual, li2017learning, chaudhry2018riemannian, aljundi2018memory} avoid updating weights that are important for predictions of previous tasks. Rehearsal \citep{rebuffi2017icarl, wu2019large, douillard2020podnet} methods keep a small portion of the previous data and replay it during training for the new task to prevent forgetting. Currently, rehearsal approaches perform better than regularization-based ones because they keep examples of previous tasks in memory. However, both types of algorithms do not completely prevent forgetting \citep{masana2020class}, even if the task-ID is given during inference. Therefore, the authors recently developed an architectural method called Continual-Prune-and-Select (CP\&S) \citep{dekhovich2022continual} that exhibits no forgetting when the task-ID is given, and that performs better than state-of-the-art algorithms for different standard computer vision datasets when the task-ID is unknown during inference. Similarly to other architectural methods such as PackNet \citep{mallya2018packnet}, Piggyback \citep{mallya2018piggyback}, CLNP \citep{golkar2019continual} and SupSup \citep{wortsman2020supermasks}, CP\&S \citep{dekhovich2022continual} finds a subnetwork for every task within an original network and uses only one subnetwork for each task. The task-related subnetwork is found by a novel iterative pruning strategy \citep{dekhovich2021neural}, although different methods could be considered \citep{ramanujan2020s, han2015learning, hu2016network}.

Here the focus is to be the first to apply continual learning strategies to mechanics problems, in particular for plasticity modeling. These cases are expected to involve known task-IDs, as illustrated by the previously invoked example of Team A and B that were aiming to model the behavior of two different materials because each team already knows \textit{a priori} the material of interest to them (known task-ID). Therefore, this article only needs to consider the simpler formulation of our CP\&S method.

Unlike the standard computer vision problems for which CP\&S was originally implemented and compared \citep{dekhovich2022continual}, here we apply CP\&S to model irreversible material behavior. Without loss of generality, we focus on plasticity simulations but other history- and time-dependent phenomena such as damage or visco-elasticity could be considered. In fact, CP\&S can be applied to other classification or regression tasks (they do not need to be tasks in Mechanics). Given the time and/or history-dependency of these problems, considering neural network architectures with recurrent units facilitates the learning process.

Recurrent neural networks (RNNs) \citep{rumelhart1985learning, jordan1997serial} typically deal with sequential data, e.g. in natural language processing problems \citep{mikolov2010recurrent} or voice recognition \citep{sak2014long}. However, they suffer from vanishing and exploding gradients issues \citep{hochreiter1991untersuchungen, bengio1994learning}. Further improvements of RNNs, such as the Long Short-Term Memory (LSTM) network \citep{hochreiter1997long} and Gated Recurrent Unit (GRU) \citep{cho2014learning} solved this problem, enabling their application in different contexts. In particular, LSTMs and GRUs have been shown to learn history-dependent phenomena in mechanics \citep{mozaffar2019, chen2021b}. In this work, we use GRUs as described in our past work \citep{mozaffar2019} due to their simplicity and effectiveness when compared to LSTM, although other neural network architectures with recurrent units could be considered. Computations that occur in the GRU can be described as follows:

\begin{align}
    z_t &= \text{sigmoid} \big(W_z x_t + U_z h_{t-1} + b_z \big),\\
    r_t &= \text{sigmoid} \big(W_r x_t + U_r h_{t-1} + b_r \big),\\
    \hat{h}_t &= \tanh \big(W_h x_t + U_h (r_t \odot h_{t-1}) + b_h \big),\\
    h_t &= z_t \odot  \hat{h}_t + (1 - z_t) \odot h_{t-1},
\end{align}

\noindent where $h_0 = 0$, \  $x_t$ is an input vector, $t = 0, 1, \ldots, T$. Matrices $W_z$, $W_r$, $W_h$, $U_z$, $U_r$, $U_h$ are learnable weights and $b_z, b_r, b_h$ are learnable biases. In this work, the PyTorch \citep{paszke2019pytorch} implementation of GRUs is used to conduct the experiments consisting of simplified representative volume elements of materials with different microstructures, as described in Section \ref{sec:plates}.

Independently of the architecture chosen for the neural network, we advocate here that architectural continual learning algorithms such as CP\&S \citep{dekhovich2022continual} enable Cooperative Data-Driven Modeling (CDDM). The method creates different subnetworks within an artificial neural network that are associated to particular tasks without forgetting past tasks. For example, in our context each subnetwork is associated to a specific material model for a class of materials with a given microstructure, constituent properties and external conditions (see e.g. \cite{bessa2017}). Architectural methods, however, are limited by the availability of free neural connections that can be trained for a new task. A particular advantage of CP\&S is that it is based on the NNrelief iterative pruning strategy \citep{dekhovich2021neural} which was developed aiming to create sparser subnetworks (using fewer connections) than other pruning methods such as magnitude pruning \citep{han2015learning} or neural pruning \citep{hu2016network}. Additionally, although CP\&S was implemented in the original article for convolutional and fully connected networks, it can easily be adapted to GRU networks as shown in this work. 

\begin{figure}[ht]
    \centering
    \includegraphics[width=0.95\textwidth]{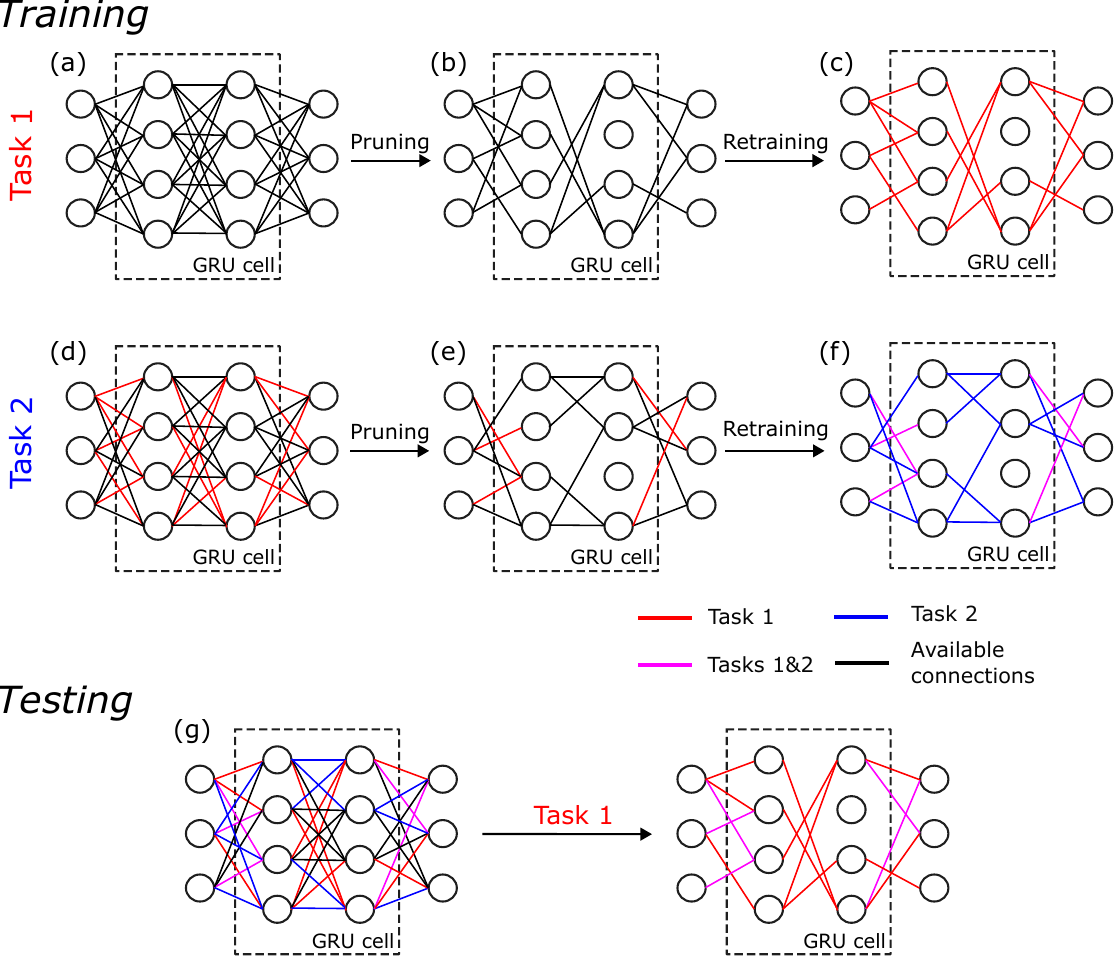}
    \caption{Overview of the proposed CDDM approach.}
    \label{fig:cddm_pipeline}
\end{figure}

CP\&S is based on a set of simple steps to create a group of overlapping subnetworks, each of them learning a particular task without disrupting the knowledge accumulated by the other subnetworks. Importantly, the subnewtorks can (and usually do) share knowledge among them by sharing connections that are useful to each other. This mechanism allows to learn different tasks, and transfer knowledge between them but avoids forgetting, unlike transfer learning methods.

Overall, cooperative data-driven modeling via the proposed CP\&S method is described as follows:

\begin{itemize}
    \item Training:
    \begin{enumerate}
        \item At the beginning of the learning process, initialize the entire neural network with random weights (see Figure \ref{fig:cddm_pipeline}a).
        \item Set the hyperparameters (architecture, optimizer, learning rate, weight decay, pruning parameters, etc.).
        \item \label{CPS_NNrelief} For a task $T_i$ (e.g. learning a new material), create a subnetwork $\mathcal{N}_i$ that is associated to that task:
        \begin{enumerate}
            \item \label{CPS:prune_step} Use the NNrelief \citep{dekhovich2021neural} algorithm to prune connections from the entire neural network (see Figure \ref{fig:cddm_pipeline}b,  details in Appendix \ref{appendix:nnrelief}).
            \begin{itemize}
                \item Remark 1: Every connection can be pruned, whether they are already part of another subnetwork or not.
                \item Remark 2: Pruning is controlled by hyperparameter $\alpha \in (0, 1)$ that represents the amount of information that is maintained coming out of the neurons after removing the connections. So, $\alpha=0.95$ indicates that 95\% of the signal coming out of the neurons is kept on average for the training set of task $T_i$ after removing the connections. The amount of information is measured according to a metric called \textit{importance score} -- see the original article \citep{dekhovich2021neural} or \ref{appendix:nnrelief} for details.
                \item Remark 3: Allowing to prune any connection of the entire network (whether belonging to a previous subnetwork or not) is crucial because it provides a mechanism to keep connections that are useful for performing the current task but remove the ones that are not. If some connections are kept from other subnetworks, then there is potential for knowledge transfer.
            \end{itemize}
            \item Among the remaining connections, i.e. the ones that are not pruned, retrain the ``free'' connections but do not update the ones that are part of other subnetworks (see Figure \ref{fig:cddm_pipeline}c).
            \begin{itemize}
                \item Remark 4: By refraining from updating connections that are part of other subnetworks, the CP\&S algorithm avoids forgetting past tasks because the other subnetworks are not affected by the current training of this subnetwork. Instead, only connections that have not been used previously by any subnetwork are updated so that new knowledge can be accumulated, as needed to solve a new task.
            \end{itemize}
            \item Assess performance of this subnetwork ($\mathcal{N}_i$) on the dataset for this task ($T_i$). If error is not acceptable, go to \ref{CPS:prune_step} and iterate again until a maximum number of iterations $n$ is met. Else save the connections that form this new subnetwork into mask $\mathbf{M}_i$. The weights belonging to this subnetwork will no longer be updated when training another subnetwork later.
        \end{enumerate}
        \item Consider the next task $T_{i+1}$ by finding a new subnetwork ($\mathcal{N}_{i+1}$) as described in step \ref{CPS_NNrelief} (see Figure \ref{fig:cddm_pipeline}d-f)).
    \end{enumerate}
    \item Inference (testing):
    \begin{enumerate}
        \item For the given test point $x$ from task $T_i$, select subnetwork $\mathcal{N}_i$ (see Figure \ref{fig:cddm_pipeline}g).
        \begin{itemize}
            \item Remark 5: in this paper, task selection from the data is not necessary, because practitioners already know what is the material model that they want to consider. Otherwise, see the original article for the task incremental learning scenario \citep{dekhovich2022continual}. 
        \end{itemize}
        \item Make a prediction $y = \mathcal{N}_i(x)$.
    \end{enumerate}
\end{itemize}

We highlight from the above description that the performance of CP\&S is controlled only by two hyperparameters: the number of pruning iterations $n$  and the pruning parameter $\alpha \in (0, 1)$. A low value of $\alpha$ will cause more connections to be pruned (information compression) but also lower expressivity (worst performance of the subnetwork after pruning and retraining for $n$ iterations). The subnetwork is represented by a binary mask $\mathbf{M}$, where every active connection is represented by 1 and where every inactive connection is represented by 0. The weights and biases that are first assigned to a particular task are not updated for subsequent tasks -- ensuring that the original subnetwork for which these connections were trained remains unchanged by the training process of a new subnetwork. This strategy was proven effective in the context of computer vision, and it should remain valid for computational mechanics applications.

An important disadvantage of CP\&S is that fixing the parameters associated with a trained subnetwork can quickly exhaust the number of ``free'' connections available to be trained in subsequent tasks. In other words, the artificial neural network can be ``saturated'' after several tasks have been learned. This is reported in the original investigation \cite{dekhovich2022continual}. However, we note that CP\&S does not have limitations on the type of neural network architecture to be considered. Therefore, the choice of architecture depends only on the solvable problem. In this paper, we consider the case when all tasks have the same input and output dimensions. This condition can be relaxed if one uses separate task-specific input and output heads. In this case, shared parameters are the ones that are in the intermediate layers.

\section{First case study: learning plasticity laws for different microstructures}\label{sec:plates}

The first and simpler case study considers non-periodic material domains subjected to uniform displacements at the boundary. The goal is to create constitutive models from these material domains similarly to the original publication on this subject \citep{bessa2017}, but including history-dependent behavior \cite{mozaffar2019}. We consider four different materials with different porous microstructures but the same matrix phase -- a simple von Mises plasticity model. This canonical example can be extended to more complex cases, including realistic microstructures, as shown in case study 2.

We present two case studies to demonstrate the usefulness of cooperative data-driven modeling (CDDM). The aim of CDDM is to create a model that is capable of performing multiple tasks but that also uses past knowledge from different tasks to decrease the number of training points required to learn a new task. Therefore, CDDM should achieve better performance than conventional training (without cooperation) for the same number of training points. In addition, the fact that CDDM uses the same neural network should make it more efficient in terms of the number of new parameters needed to learn a new task.

\subsection{Data generation}\label{subsec:data}
\begin{figure}[ht]
  \centering
  \includegraphics[scale=0.3]{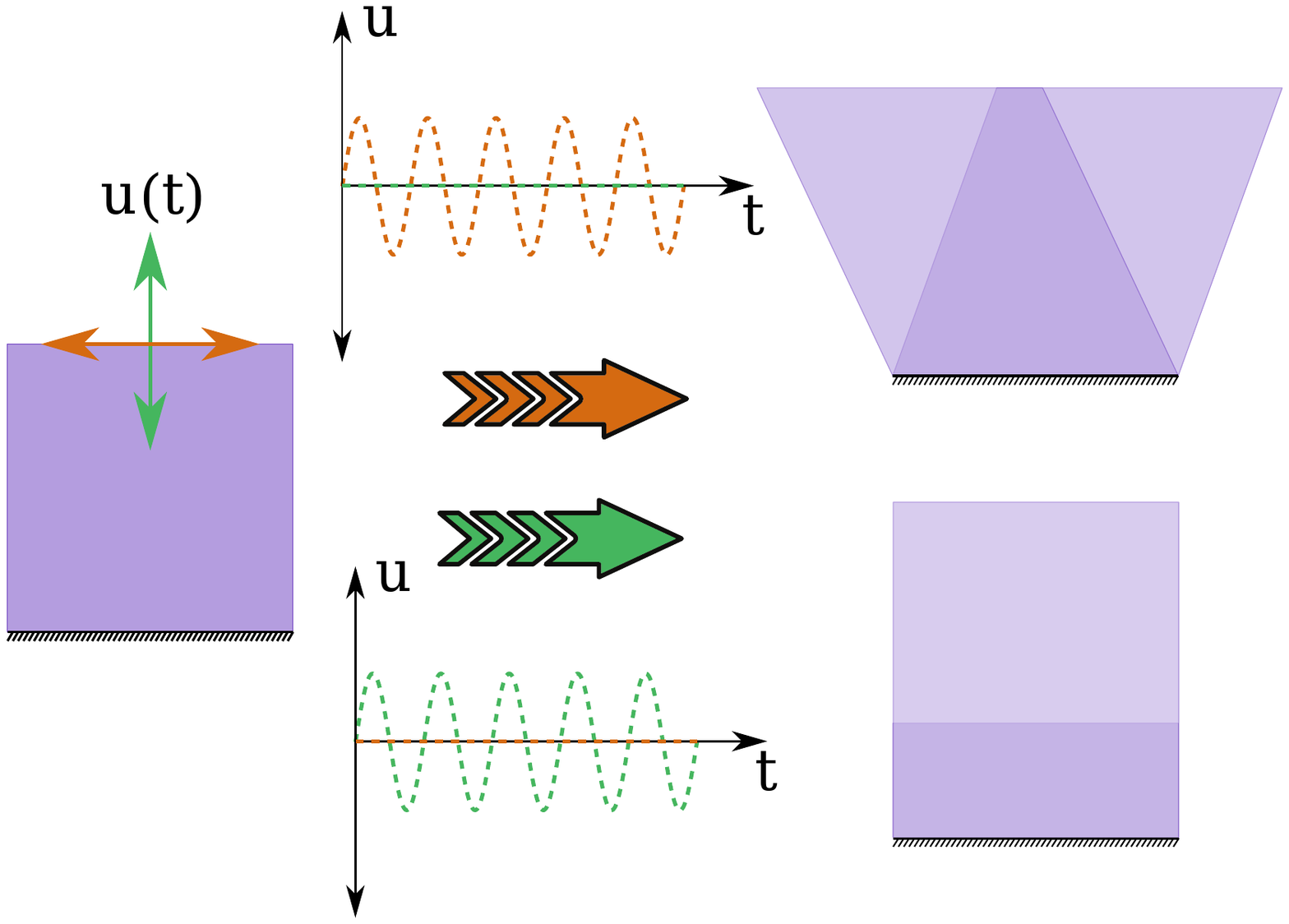}
  \caption{A square domain fixed on the bottom and displaced on the top. Displacement is done in a pseudo-time.}
  \label{fig:problem}
\end{figure}

\begin{figure}[ht]
  \centering
  \begin{subfigure}[b]{0.22\textwidth}
    \centering
    \includegraphics[width=\textwidth]{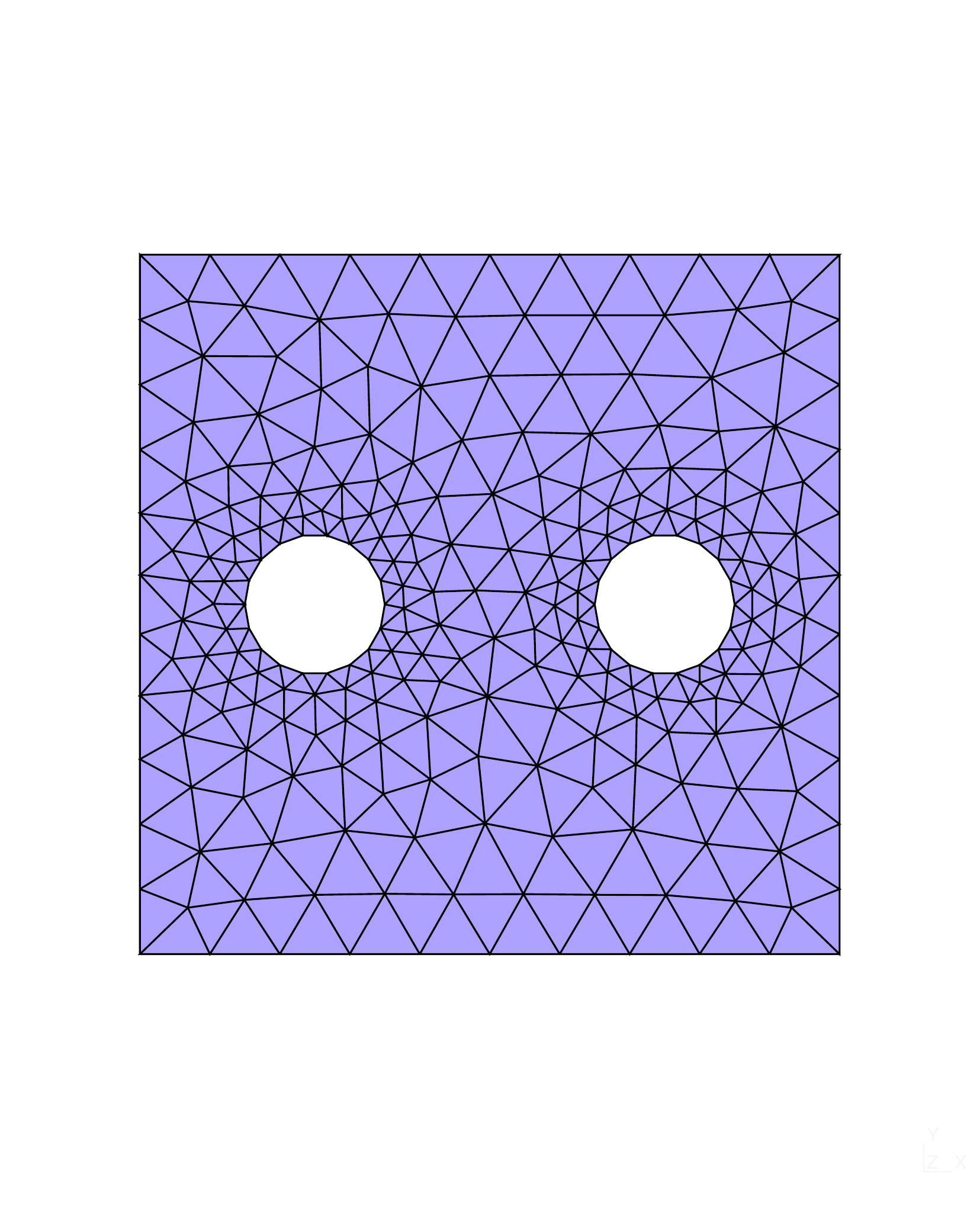}
    \caption{Task A.}
    \label{fig:task1}
  \end{subfigure}
  \hfill
  \begin{subfigure}[b]{0.22\textwidth}
    \centering
    \includegraphics[width=\textwidth]{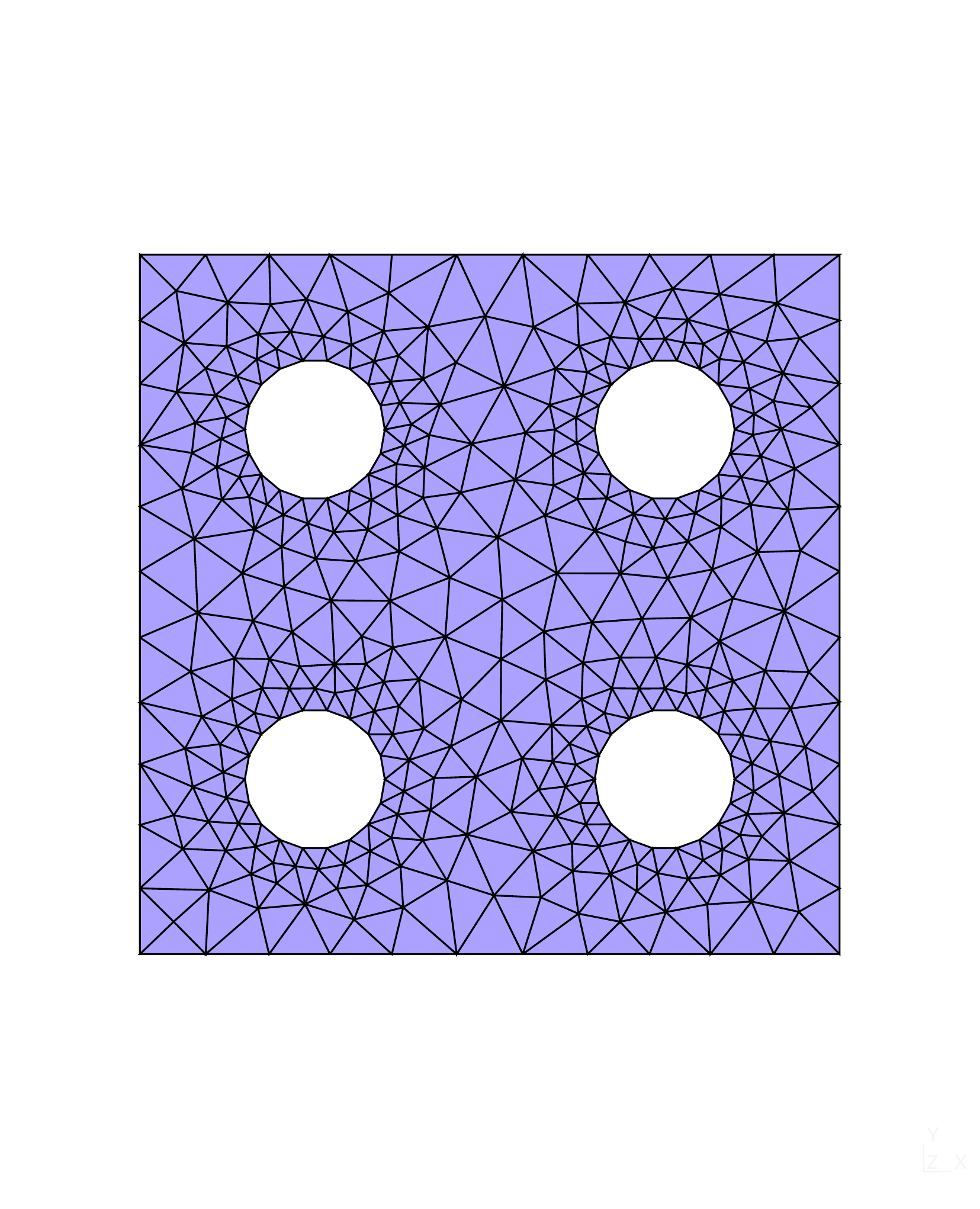}
    \caption{Task B.}
    \label{fig:task2}
  \end{subfigure}
  \hfill
  \begin{subfigure}[b]{0.22\textwidth}
    \centering
    \includegraphics[width=\textwidth]{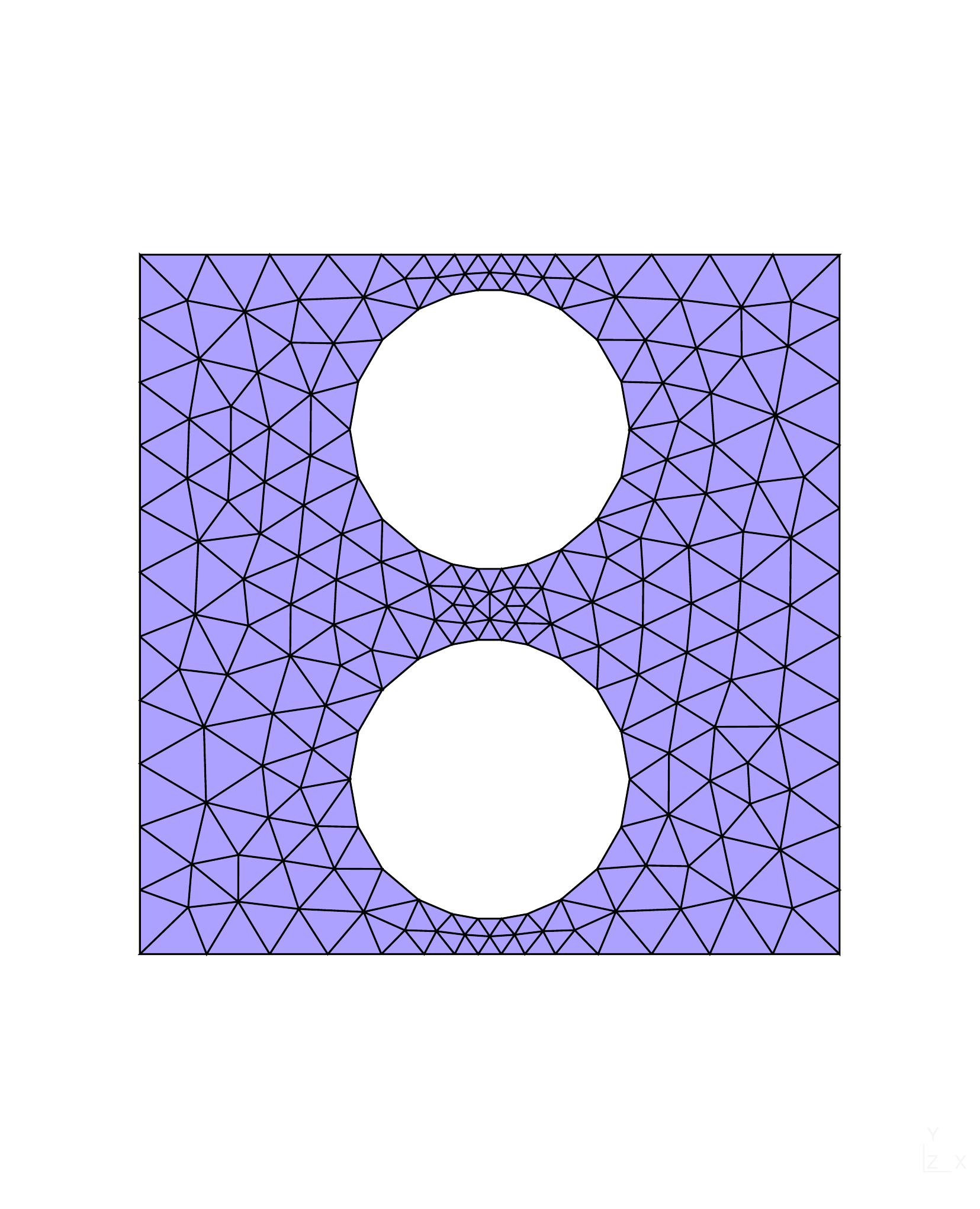}
    \caption{Task C.}
    \label{fig:task3}
  \end{subfigure}
  \hfill
  \begin{subfigure}[b]{0.22\textwidth}
    \centering
    \includegraphics[width=\textwidth]{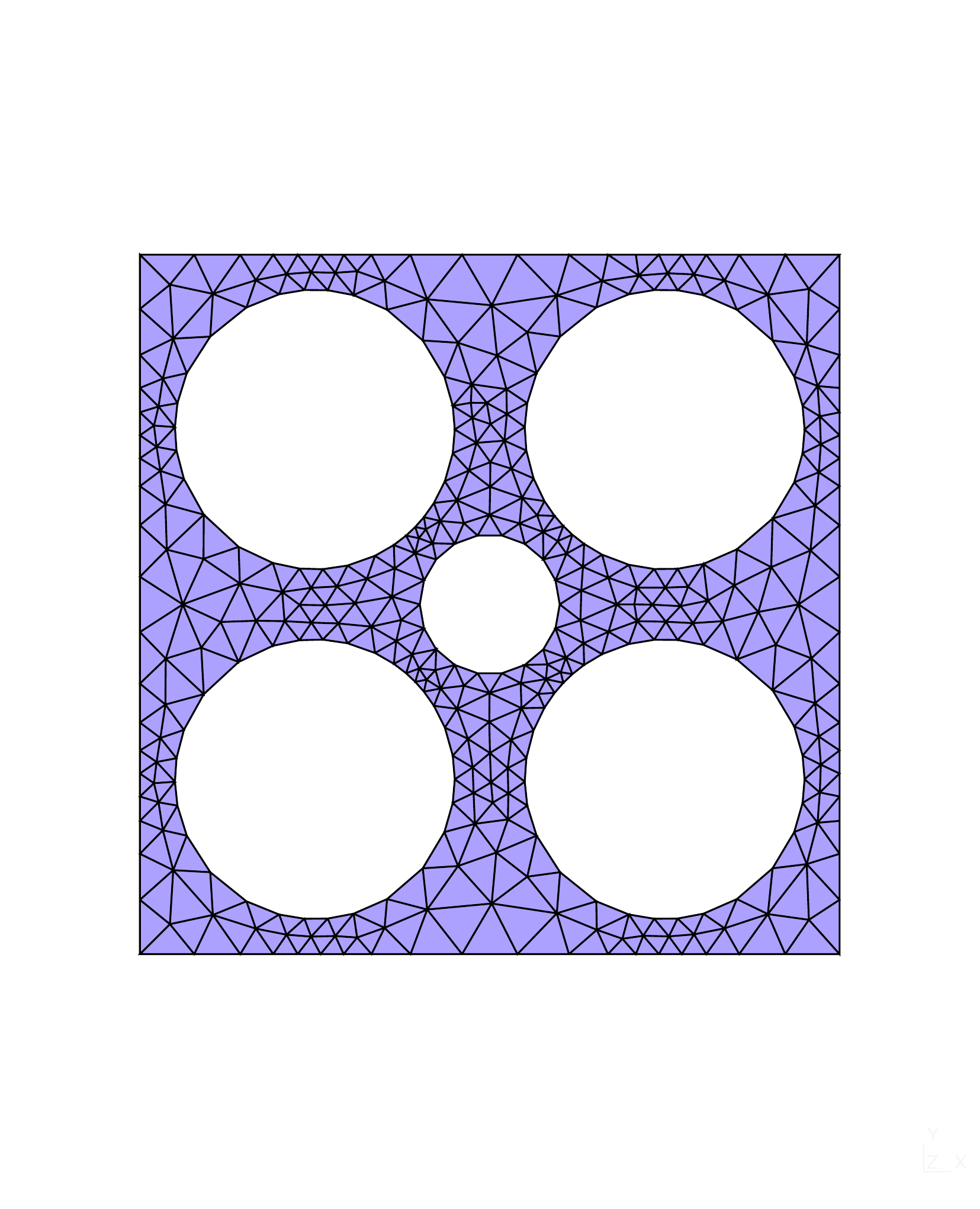}
    \caption{Task D.}
    \label{fig:task4}
  \end{subfigure}
  \caption{Different domains introduced as different tasks.}
  \label{fig:tasks}
\end{figure}

% Problem explanation
An elastoplastic von Mises material with hardening is investigated to create a path-dependent problem. A fixed-sized square is utilized as a domain and holes of varying sizes and locations are placed inside the domain to create different tasks (see Figure \ref{fig:problem}). The tasks originating from different domains can be seen in Figure \ref{fig:tasks}. For all the tasks the bottom part of the domain is fixed and the top part is deformed according to a uniform displacement in $\dispx$ and $\dispy$. We reinforce that Case Study 1 considers simpler finite element analysis because we want to facilitate the dissemination of the methodology and use fully open-source software -- facilitating replication of the example and easy adaptation to other scenarios.

After applying the boundary conditions, the average Cauchy strain ($\boldsymbol{\avgstrain}$) and Cauchy stress ($\boldsymbol{\avgstress}$) measures are obtained for the deformation path. Then the learning problem for a single task can be defined as,
\begin{equation}
  \boldsymbol{\avgstress} = f(\boldsymbol{\avgstrain}),
\end{equation}
where a machine learning model is utilized to find the relationship $f: \boldsymbol{\avgstrain} \mapsto \boldsymbol{\avgstress}$, which is a supervised regression problem with the history of the average strain components as an input and the average stress components as an output.

% Data generation

\begin{figure}[ht]
    \centering
    \includegraphics[width=\textwidth]{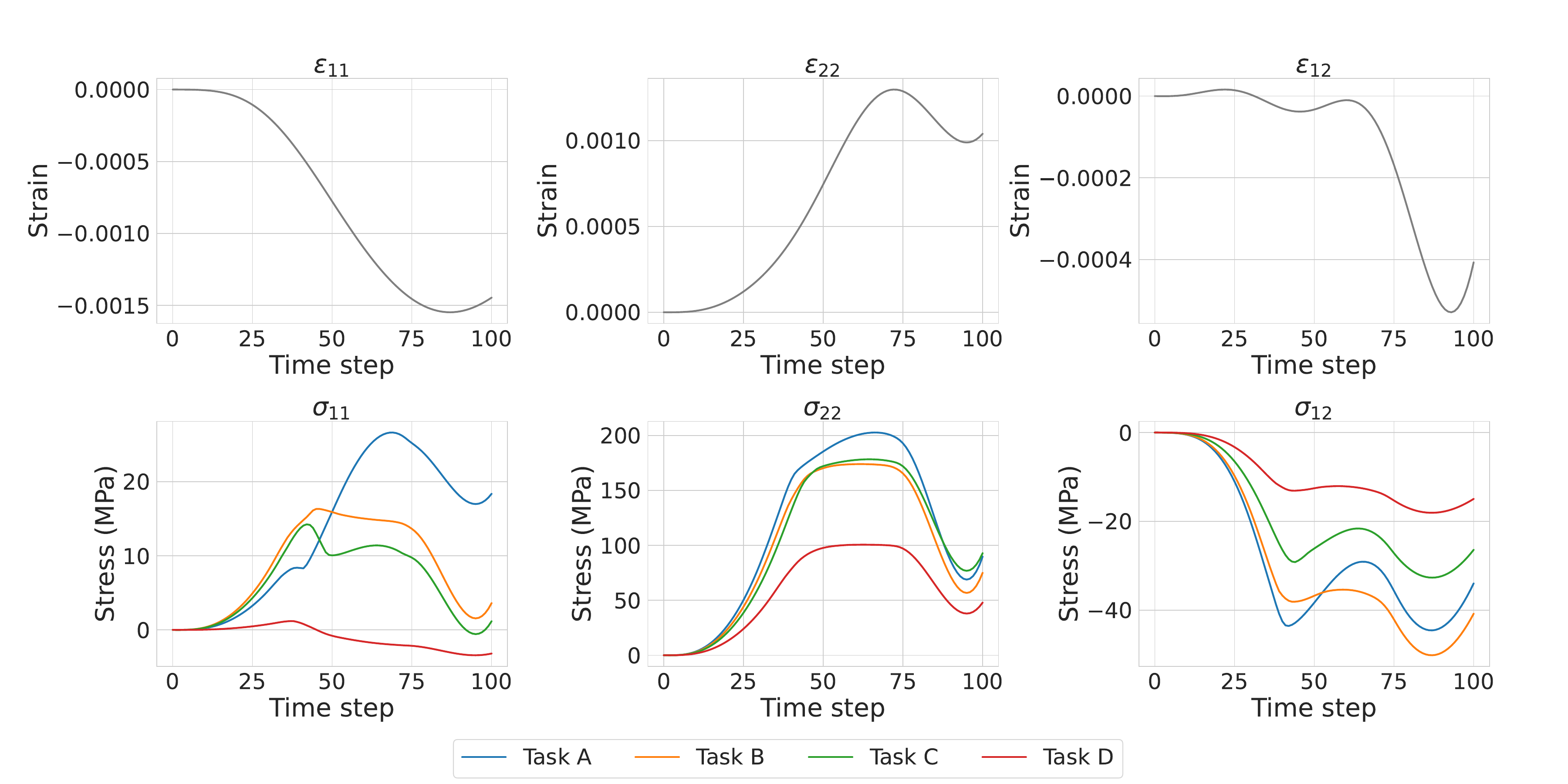}
    \caption{An example of strain and stress paths.}
    \label{fig:plates_paths}
\end{figure}

The four material domains are subjected to the same 1000 paths of deformation. These 1000 paths are obtained from 100 end displacement values of the top boundary of the domain that are sampled from a Gaussian Process posterior that is conditioned on 20 displacement values sampled from a uniform distribution for each path. Then the average stress obtained for each path is calculated via the finite element method -- see for example Figure \ref{fig:plates_paths}. The dataset is made available, and the simulations used to obtain it can be replicated via the source code. We highlight that different strategies for sampling the paths of plasticity material laws have been proposed \citep{wu2020recurrent} and that this choice can affect the number of paths needed to train the neural network up to the desired accuracy. Nevertheless, each presented task was subjected to the same deformation paths to calculate the domain-specific average stress and strain values. The data is generated using FEniCS \citep{logg2012a}.

\begin{figure}[ht]
    \centering
    \includegraphics[width=\textwidth]{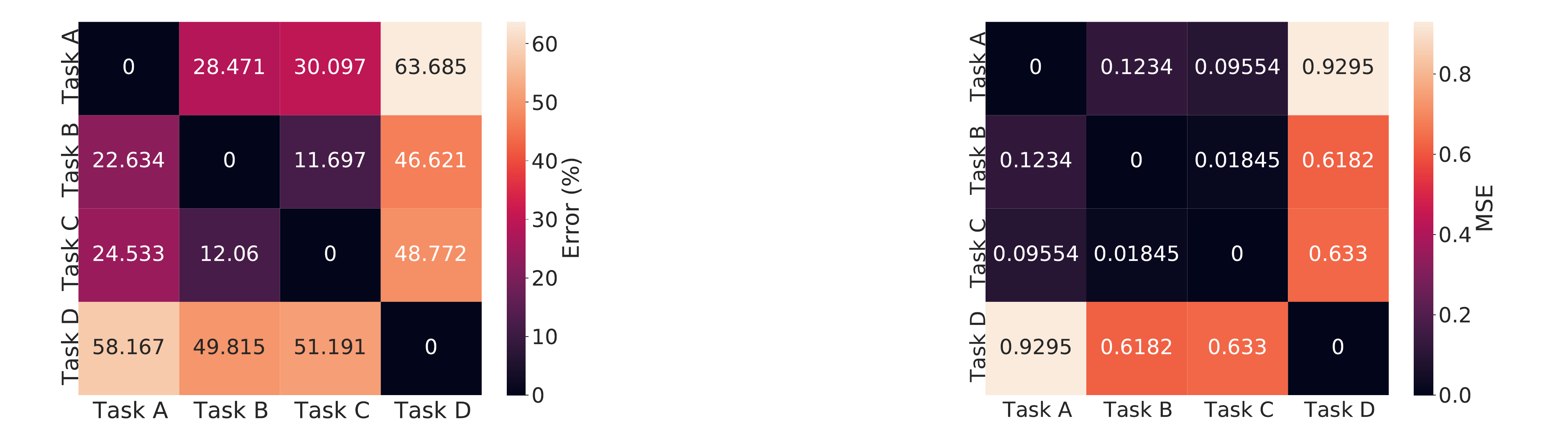}
    \caption{The difference between the stress train data for considered tasks, calculated using the error measurement (left) and the MSE measurement (right).}
    \label{fig:data_comparison}
\end{figure}

In order to illustrate the difference in stresses obtained for the four tasks, Figure \ref{fig:data_comparison} shows the error according to Eq. \ref{eqn:error} and the mean-squared-error (MSE) between stresses for the training data for all tasks. We scale all the data since the model receives it in a scaled format. Here we use standard scaling where we remove the mean from every input feature and divide it by the standard deviation of the training set. The figure clarifies that Tasks A, B and C are more similar to each other than Task D.

\subsection{Experiments setup} \label{sec:exp_setup}

% Continual Learning Problem (I think we should here define the overall continual learning problem)

To evaluate the performance of CDDM we measure the error $E_i$ on every test path as follows:

\begin{equation}\label{eqn:error}
    E_i = \frac{1}{3} \Biggl( \frac{\lvert\lvert\sigma^i_{11} - \hat{\sigma}^i_{11} \rvert\rvert_2}{\lvert\lvert\sigma^i_{11}\rvert\rvert_2} + \frac{\lvert\lvert\sigma^i_{22} - \hat{\sigma}^i_{22}\rvert\rvert_2}{\lvert\lvert \sigma^i_{22} \rvert\rvert_2} + \frac{\lvert\lvert\sigma^i_{12} - \hat{\sigma}^i_{12}\rvert\rvert_2}{\lvert\lvert\sigma^i_{12}\rvert\rvert_2} \Biggl) \cdot 100\%,\ i=1, 2, \ldots, N,
\end{equation}
where $N$ is the number of testing points, $\lvert\lvert \cdot \rvert\rvert_2$ is a L2 norm, $\hat{\sigma}^i_{11}, \hat{\sigma}^i_{22}, \hat{\sigma}^i_{12} \in \mathbb{R}^t$ are predicted stress components, and $\sigma^i_{11}, \sigma^i_{22}, \sigma^i_{12} \in \mathbb{R}^t$ are the test ones.

Then, we compute the average over all $N$ test points to compute the final test error:
\begin{equation}
    Err = \frac{1}{N} \sum_{i=1}^N E_i.
\end{equation}

The hyperparameters that we use to train a neural network are shown in Table \ref{tb:haperparameters}. We use Adam \citep{kingma2014adam} optimizer to train the model with the mean-squared-error (MSE) loss function. 

\begin{table}
    \centering
    \caption{Training hyperparameters.}\label{tb:haperparameters}
    \begin{tabular}{p{0.7in} p{0.7in}  p{0.7in} p{0.7in} p{0.8in} p{0.7in}}
    \toprule
     training epochs & learning rate & weight decay & pruning iterations & pruning parameter $\alpha$ & retraining epochs  \\ % Table header row
    \midrule
      1000 & 0.01 & $10^{-6}$ & 1 & 0.95 & 200 \\
    \bottomrule
    \end{tabular}
\end{table}

As it is common in continual learning literature \citep{masana2020orderings}, we test the approach with different task orderings. Overall, we consider four orderings for the case of four tasks in a sequence:

\begin{itemize}
  \item ordering 1: Task A $\rightarrow$ Task B $\rightarrow$ Task C $\rightarrow$ Task D;
  \item ordering 2: Task B $\rightarrow$ Task D $\rightarrow$ Task A $\rightarrow$ Task C;
  \item ordering 3: Task C $\rightarrow$ Task A $\rightarrow$ Task D $\rightarrow$ Task B;
  \item ordering 4: Task D $\rightarrow$ Task C $\rightarrow$ Task B $\rightarrow$ Task A.
\end{itemize}

\subsection{Results}\label{subsec:results}

\begin{figure}[ht!]
  \centering
  \begin{subfigure}[b]{\textwidth}
    \centering
    \includegraphics[width=\textwidth]{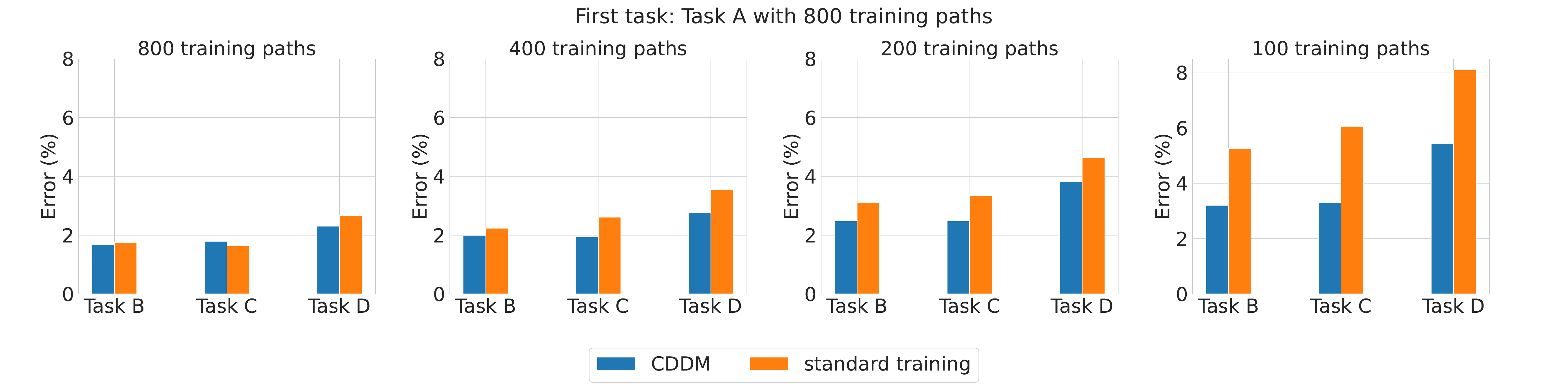}
	\caption{Ordering 1.}\label{fig:4tasks_order1}
  \end{subfigure}%
  \vspace{1em}
  \begin{subfigure}[b]{\textwidth}
    \centering
    \includegraphics[width=\textwidth]{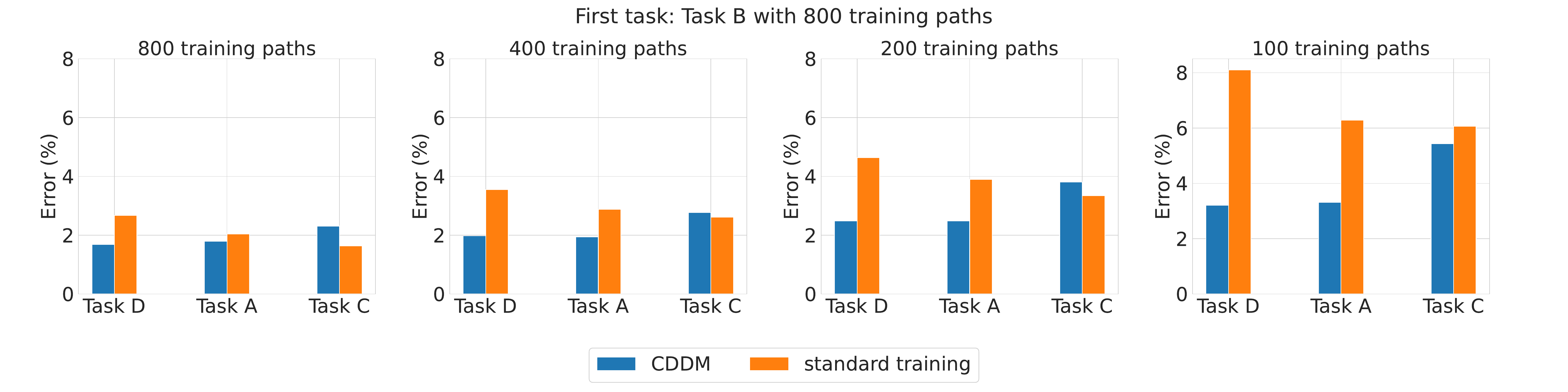}
	\caption{Ordering 2.}\label{fig:4tasks_order2}
  \end{subfigure}%
  \vspace{1em}
  \begin{subfigure}[b]{\textwidth}
    \centering
    \includegraphics[width=\textwidth]{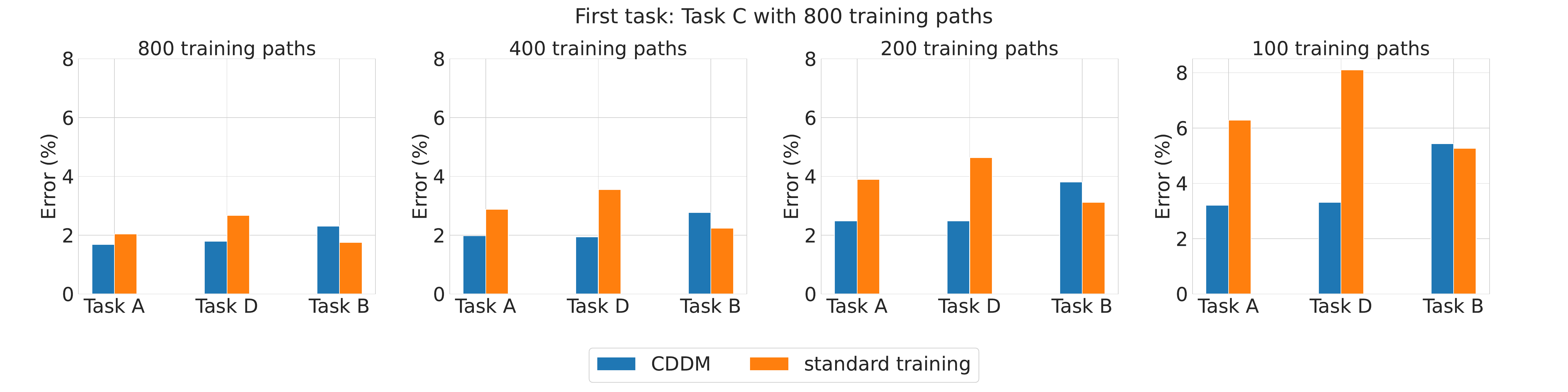}
	\caption{Ordering 3.}\label{fig:4tasks_order3}
  \end{subfigure}%
  \vspace{1em}
  \begin{subfigure}[b]{\textwidth}
    \centering
    \includegraphics[width=\textwidth]{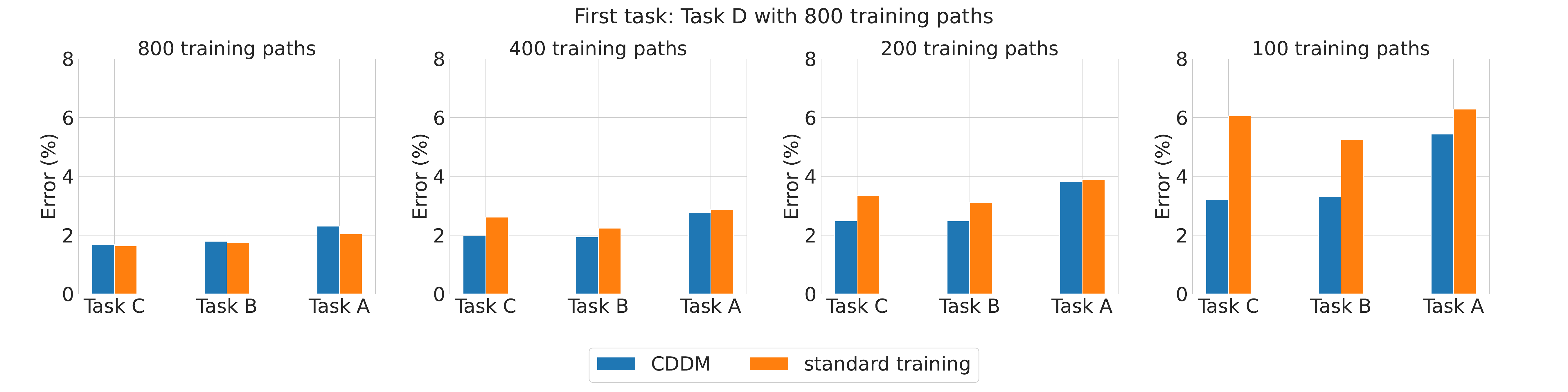}
	\caption{Ordering 4.}\label{fig:4tasks_order4}
  \end{subfigure}
  \caption{First case study: CDDM results on orderings 1-4.}
  \label{fig:4orderings}
\end{figure}

Firstly, we train the GRU with 2 cells and a hidden size of 128 in the sequence of four tasks. We train the first task with 800 training paths, and for each of the following tasks, we consider the cases of 800, 400, 200 and 100 training paths. We compare these results with the conventional case (non-cooperative) where every new task is trained with the same GRU but independently of the other tasks. Figure \ref{fig:4orderings} shows this comparison, where the blue bars refer to the cooperative model (CDDM) and the orange ones to the conventional case (standard or non-cooperative training). The test error is computed using Eq. \ref{eqn:error}. It is clear that CDDM significantly outperforms standard training when we decrease the number of training points. This effect is consistent across all four orders, independently of which task is considered to be the first. The main advantage of CDDM is that the pretrained parameters have an accumulative effect on future tasks. This multi-transfer effect has more significance under the low-data regime (e.g., 100 training paths). Also, the set of parameters depends on the order in which the tasks are learned. In Table \ref{tab:4orderings_result}, we present the average error for every task when considering a different number of training points (number of paths); note that this error is the average over the four task orderings. We can clearly see that CDDM performs better than standard training for tasks 2 -- 4 on average.

\begin{table}[ht!]
\caption{Test error ($\%$) averaged over four considered task orders.}
\label{tab:4orderings_result}
  \centering
    \begin{tabular}{p{0.8in}cccccc}
        \toprule
           & task 1 &  & \multicolumn{4}{c}{tasks 2--4}  \\
           \cline{4-7}
           & 800 paths & & 800 paths & 400 paths & 200 paths & 100 paths\\
        \midrule 
        standard training & \textbf{2.02} & & 2.02 & 2.82 & 3.75 & 6.43 \\
        \midrule 
        CDDM & 2.14 & & \textbf{1.92} & \textbf{2.18} & \textbf{2.84} & \textbf{3.97}\\
        \bottomrule
    \end{tabular}
\end{table}

Moreover, we learn all four tasks with one network, while four separate networks are necessary for the conventional case. Therefore, using fewer parameters and achieving better performance. As clarified next, this is explained by the knowledge transfer that happens between subnetworks. Also, it should be noted that the network still has free space to learn future tasks, although saturation would occur soon if more tasks were considered because the neural network is small.

\begin{figure}[ht]
    \centering
    \includegraphics[width=0.95\textwidth]{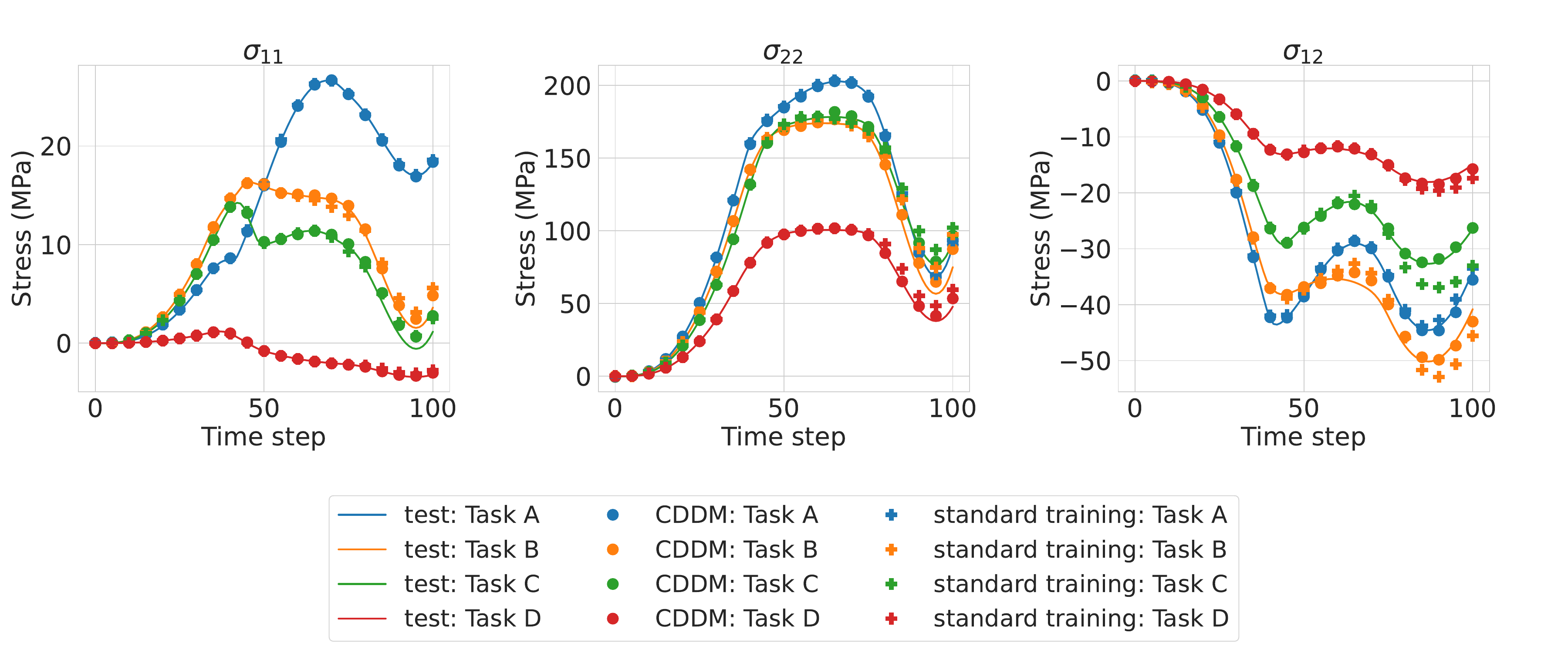}
    \caption{First case study: 800 training paths for task 1 and 200 training paths for tasks 2-4.}\label{fig:4tasks_order4_time-stress_800-200}
\end{figure}

In Figure \ref{fig:4tasks_order4_time-stress_800-200}, we show the prediction of CDDM with 200 training paths for tasks 2-4 respectively; the first task is learned with 800 paths. We compare CDDM with standard (non-cooperative) training with one network per task. As the figure illustrates, CDDM learns the data better than conventional training and requires just one network instead of four. Hence, Figure \ref{fig:4tasks_order4_time-stress_800-200} justifies our hypothesis of using continual learning as a possible solution to tackle the data scarcity problem in history-dependent constitutive law modeling. From this figure, it is clear that the GRU is able to learn material behavior by having 200 training paths for tasks 2-4.

Overall, we observe that the continual learning strategy allows us not only to learn four different geometries with one network but also improves test error under the limited data regime. In addition, we can see that even if we have enough data, continual learning does not worsen the results significantly compared to standard training (no more than 0.5\% difference).

\section{Second case study: RVEs with periodic boundary conditions}\label{sec:rves}

In this second case study, we apply the CDDM strategy to more realistic representative volume elements (RVEs) subjected to periodic boundary conditions, as thoroughly discussed in a past work \cite{bessa2017}. The data is generated with a commercial finite element solver, but both the code to generate the data as well as the dataset itself are made available. The RVE average stress--strain response is influenced by the microstructure, properties of each material phase, and loading conditions (average strain path that is converted into a periodic displacement at the boundary). Four different RVEs are considered in order to create four tasks. Each task pertains to learning the homogenized plasticity constitutive behavior (average stress--strain response) of a corresponding RVE. These RVEs were created such that they share some similarities and significant differences by considering different microstructures and material properties while applying the same average deformation paths to generate the dataset.

\subsection{RVEs simulation}

\begin{figure}[ht]
    \centering
    \includegraphics[width=0.9\textwidth]{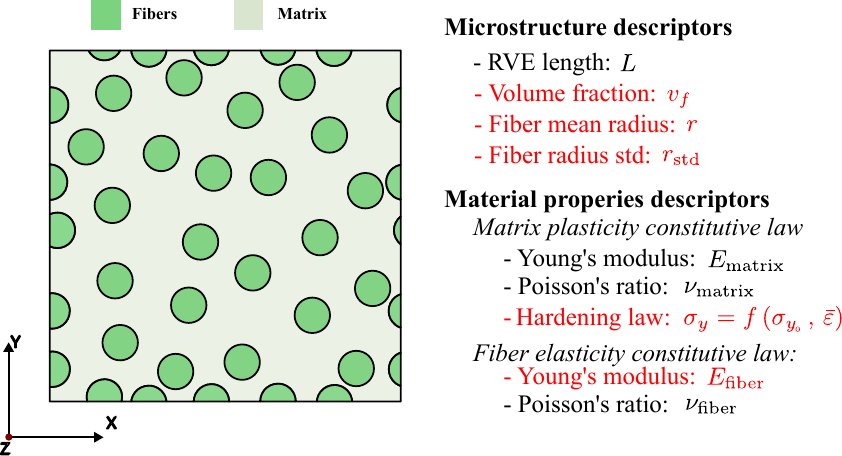}
    \caption{Illustration of commonalities and discrepancies on setting up tasks}
    \label{fig:rve_illustration}
\end{figure}

\begin{figure}[ht]
    \centering
    \includegraphics[width=\textwidth]{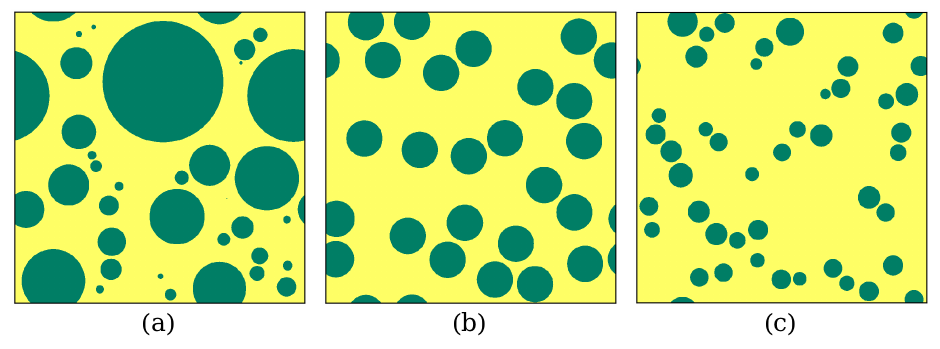}
    \caption{Schematics of different microstructure configurations: (a) $v_{f}=0.45$, $r=0.01$, $r_{\mathrm{std}}=0.003$; (2) $v_{f}=0.30$, $r=0.003$, $r_{\mathrm{std}}=0.0$ ; (c) $v_{f}=0.15$, $r=0.0015$, $r_{\mathrm{std}}=0.0003$ }
    \label{fig:task_mircostructures}
\end{figure}

Figure \ref{fig:rve_illustration} illustrates the type of two-dimensional composite material used to create the 4 RVEs (tasks). The two-phase RVEs are defined by 4 geometric descriptors: (1) RVE size, (2) fiber\footnote{We consider plane strain conditions, so we tend to view the reinforcement phase as fibers instead of particles.} volume fraction ($v_f$), (3) fiber mean radius ($r$), and (4) fiber radius standard deviation ($r_{\mathrm{std}}$). The last two descriptors are used to create circular fibers whose radius is drawn from a Gaussian distribution with the corresponding mean and standard deviation. The descriptors of the material properties are simply the elastic properties of the fibers and matrix (Young's modulus ($E$) and Poisson's ratio ($\nu$)), and the plasticity properties of the matrix (isotropic hardening law that depends on the yield stress that completely defines the von Mises yield surface). These descriptors are all defined on the right part of Figure \ref{fig:rve_illustration}, where the red font indicates the descriptors that are changed among the 4 RVEs and where the ones in black font indicate parameters that are fixed in this investigation, without loss of generality. The 4 RVEs (i.e. Tasks) are labeled A, B, C and D and the corresponding values for the descriptors are included in Table \ref{tab:rves_tasks_configurations}. For clarity, the 3 microstructures defining the RVEs (note that two RVEs have the same microstructure) are shown in Figure \ref{fig:task_mircostructures}.

\begin{table}[ht]
    \tiny
    \centering
    \caption{Parameters configuration of different tasks (Units:SI($mm$)) }
    \begin{tabular}{cccccccccc}
      \toprule
       \multirow{2}{*}{Task}  & \multicolumn{3}{c}{Mircostructure parameters} & \multirow{2}{*}{Hardening law }  & \multirow{2}{*}{$E_{\mathrm{fiber}}$} & \multicolumn{4}{c}{Fixed parameters} \\ \cline{2-4} \cline{7-10} 
       & $v_{f}$& $r$ & $r_{\mathrm{std}}$ &  & & size &$E_{\mathrm{matrix}}$&$\nu_{\mathrm{matrix}}$  & $\nu_{\mathrm{fiber}}$  \\  
       \midrule
       A & 0.45 & 0.01 & 0.003 & $\sigma_{y} = 0.5 + 0.5\Bar{\varepsilon}$ & 10 & \multirow{4}{*}{0.048 } & \multirow{4}{*}{100} &  \multirow{4}{*}{0.30} &  \multirow{4}{*}{0.19} \\
      B & 0.30 & 0.003& 0.0 &  $\sigma_{y} = 0.5 + 0.5(\Bar{\varepsilon})^{0.4}$ & 1 &  \\  
      C & 0.15 & 0.0015 & 0.0003 &  $\sigma_{y} = 0.5(1+\Bar{\varepsilon})^{\frac{1}{0.4}}$ & 1000 & \\
      D & 0.30 & 0.003 & 0.0 &  $\sigma_{y} = 3.0 + 0.5(\Bar{\varepsilon})^{0.4}$ & 1 \\
      \bottomrule
    \end{tabular}
    \label{tab:rves_tasks_configurations}
\end{table} 

\begin{figure}[ht]
    \centering
    \includegraphics[width=\textwidth]{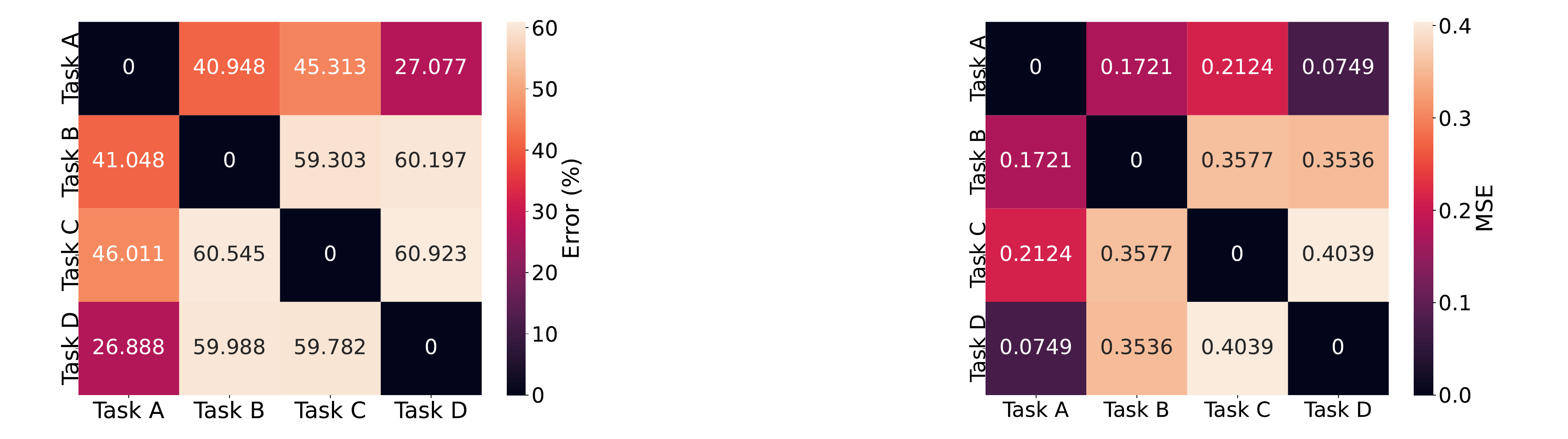}
    \caption{The difference between the stress train data for considered RVE tasks, calculated using the error measurement (left)
and the MSE measurement (right)}
    \label{fig:rve_data_comparison}
\end{figure}

As in Case Study 1, the target is to learn the average Cauchy stress ($\boldsymbol{\avgstress}$) which is dependent on the applied average strain path ($\boldsymbol{\avgstrain}$). Therefore, the learning problem can be defined as: $\boldsymbol{\avgstress} = f(\boldsymbol{\avgstrain})$, where $t$ is the pseudo-time step (load step) in the simulation defined to be 100. Meanwhile, 1000 different strain paths are generated according to a simple interpolation method. Specifically, for each average strain component ($\Bar{\varepsilon}_{11}$, $\Bar{\varepsilon}_{22}$, and $\Bar{\varepsilon}_{12}$) of a path, 8 equally spaced points are sampled within the strain path, and the quadratic interpolation method is adopted to generate the full strain path. Then, the strain path is converted into a boundary value problem of the RVE and the finite element prediction is conducted with the commercial software ABAQUS \citep{abaqus} to simulate the corresponding average stress for different tasks. In Figure \ref{fig:rve_data_comparison}, we present the difference between data computed with error metric (Eq. \ref{eqn:error}) and the mean squared error (MSE).

\subsection{Results}

\begin{figure}[ht!]
  \centering
  \begin{subfigure}[b]{\textwidth}
    \centering
    \includegraphics[width=\textwidth]{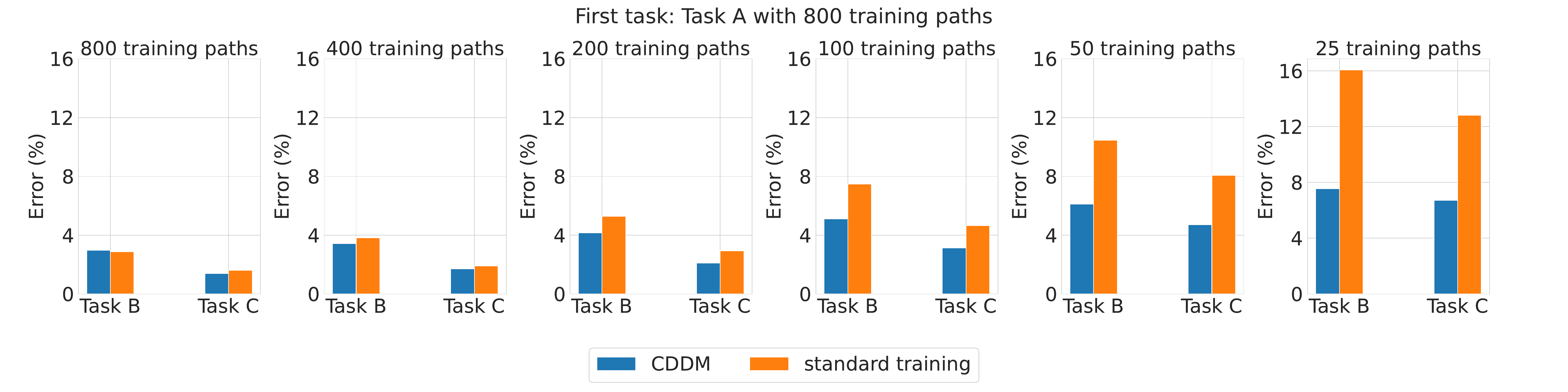}
	\caption{Ordering 1.}\label{fig:rve_3tasks_order1}
  \end{subfigure}%
  \vspace{1em}
  \begin{subfigure}[b]{\textwidth}
    \centering
    \includegraphics[width=\textwidth]{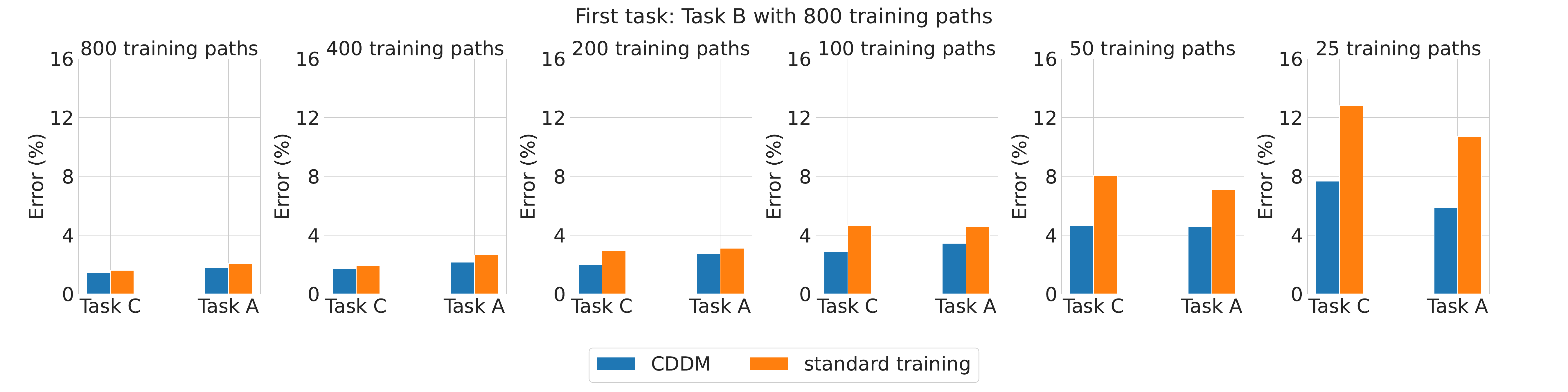}
	\caption{Ordering 2.}\label{fig:rve_3tasks_order2}
  \end{subfigure}%
  \vspace{1em}
  \begin{subfigure}[b]{\textwidth}
    \centering
    \includegraphics[width=\textwidth]{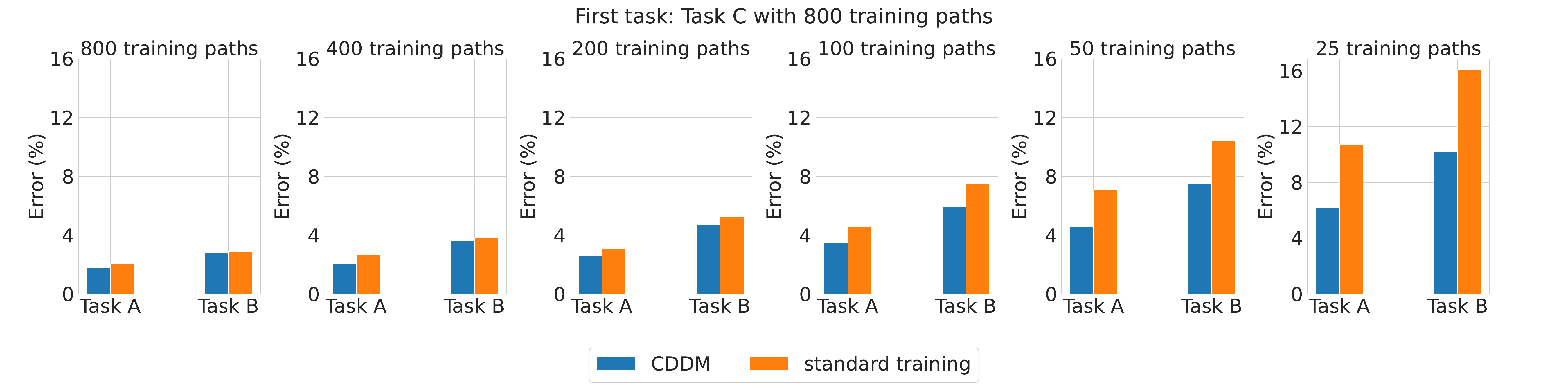}
	\caption{Ordering 3.}\label{fig:rve_3tasks_order3}
  \end{subfigure}%
  \caption{Second case study: CDDM results on orderings 1-3.}
  \label{fig:rve_3orderings}
\end{figure}

First, we train GRU on tasks A, B and C with the same hyperparameters as in Section \ref{sec:plates}. We consider three different tasks orders:
\begin{itemize}
    \item ordering 1: Task A $\to$ Task B $\to$ Task C;
    \item ordering 2: Task C $\to$ Task A $\to$ Task B;
    \item ordering 3: Task B $\to$ Task C $\to$ Task A.
\end{itemize}
In Figure \ref{fig:rve_3orderings}, we show the results for these three task orderings. It is clear that with the decrease in the number of training paths, CDDM starts to outperform conventional training. 

\begin{figure}[ht]
    \centering
    \begin{minipage}[b]{0.49\textwidth}
            \centering
            \includegraphics[width=\textwidth]{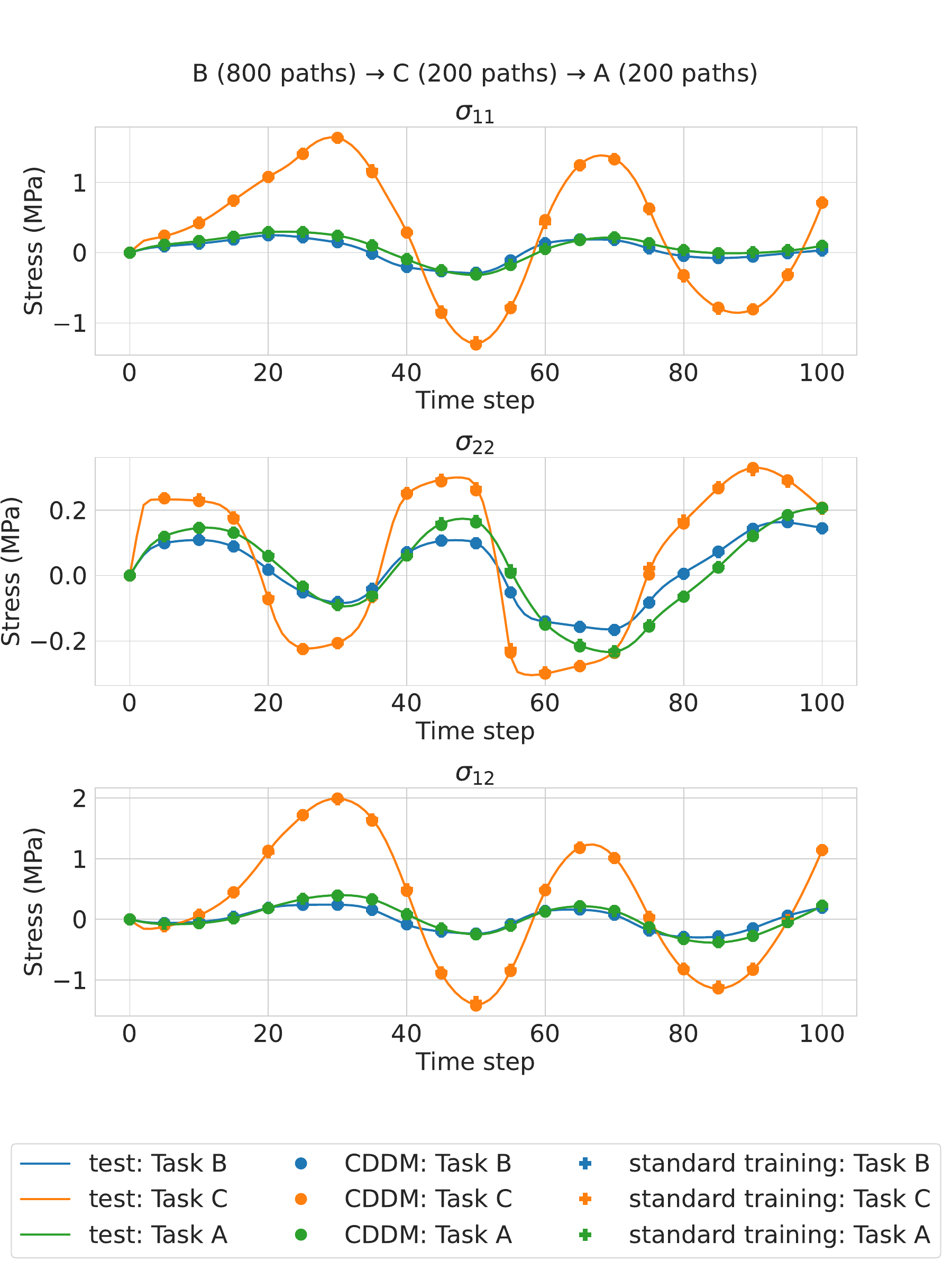}
    \end{minipage}%
    \begin{minipage}[b]{0.49\textwidth}
            \includegraphics[width=\textwidth]{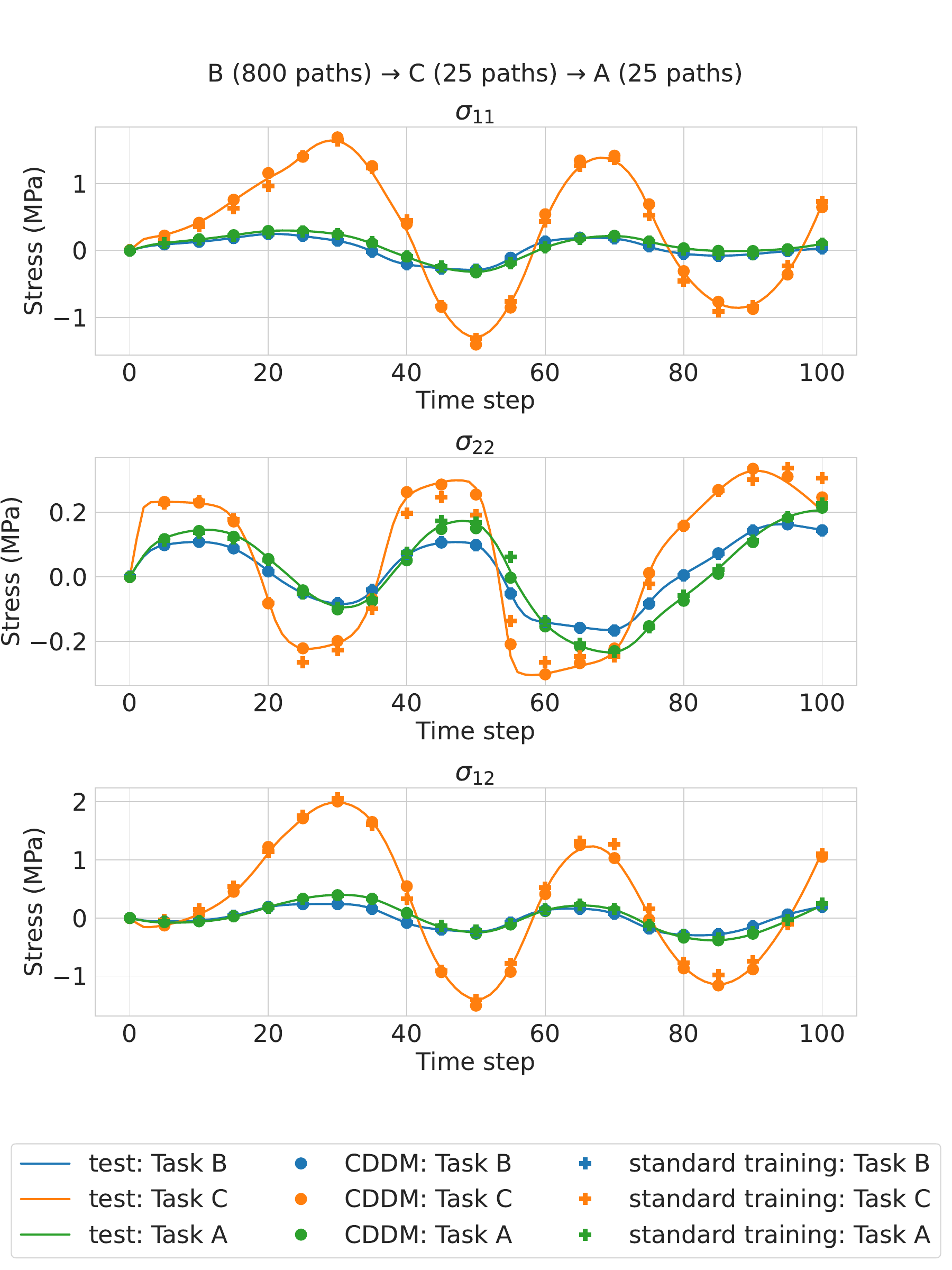}
    \end{minipage}
    \caption{Second case study: Comparison of the CDDM and standard training predictions with different numbers of training paths.}\label{fig:rve_diff_paths_predict}
\end{figure}

In Figure \ref{fig:rve_diff_paths_predict}, we present the prediction of one stress path with CDDM and standard training. The first task (task B) is trained with 800 training paths, while the next two tasks (tasks C and A) are trained using 200 training paths (\textbf{left}) and 25 training paths (\textbf{right}). If tasks C and A are learned with 200 paths, both CDDM and standard training predict stress well, however, CDDM does this with a single network. We observe that if GRU learns tasks in the cooperative approach, the prediction is more accurate than with conventional training if 25 training paths are given for the second and third tasks.

\begin{figure}[ht]
  \centering
  \begin{subfigure}[b]{\textwidth}
    \centering
    \includegraphics[width=\textwidth]{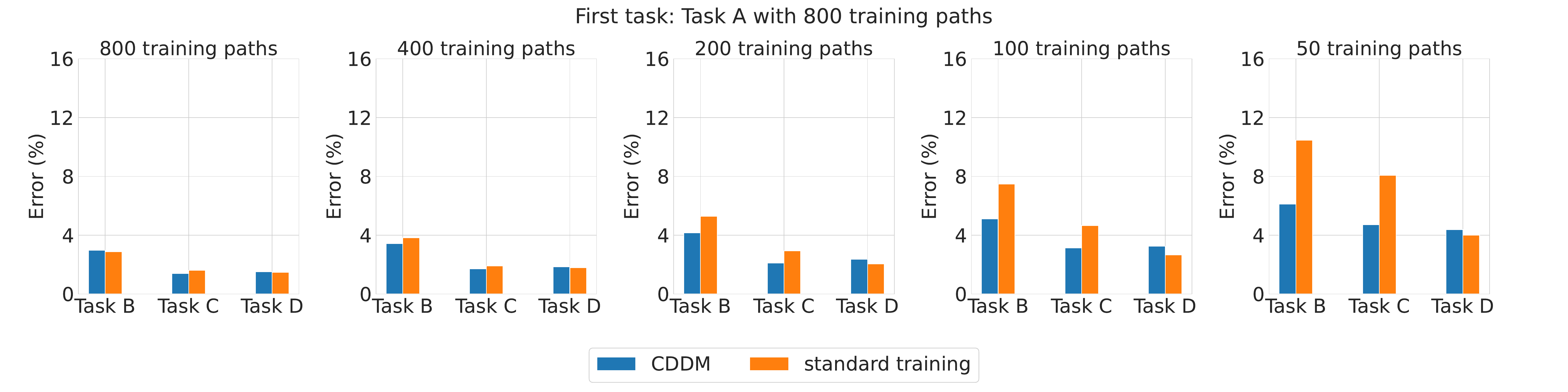}
	\caption{Task A $\to$ Task B $\to$ Task C $\to$ Task D.}\label{fig:rve_4tasks_order1}
  \end{subfigure}%
  \vspace{1em}
  \begin{subfigure}[b]{\textwidth}
    \centering
    \includegraphics[width=\textwidth]{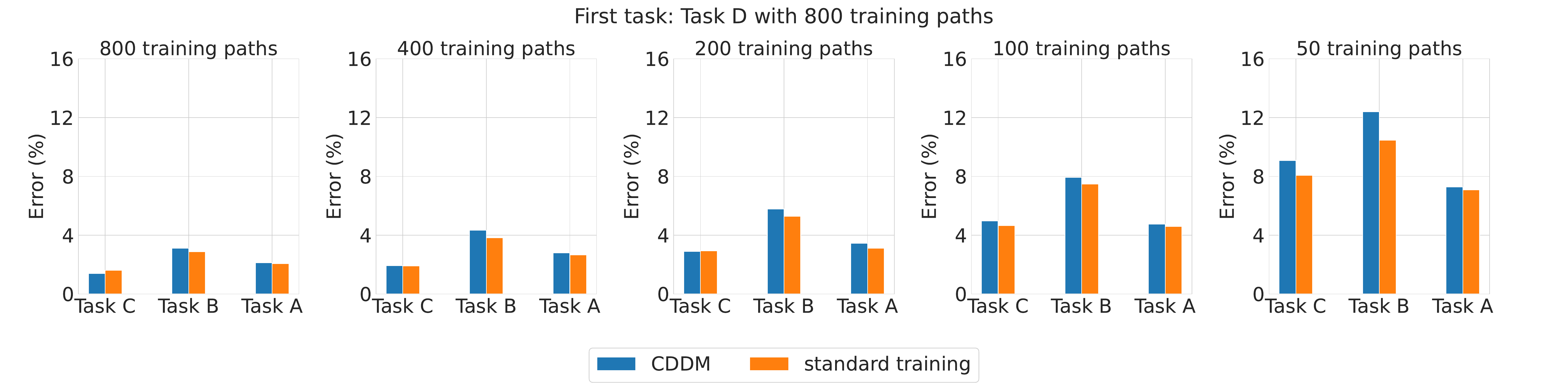}
	\caption{Task D $\to$ Task C $\to$ Task B $\to$ Task A.}\label{fig:rve_4tasks_order2}
  \end{subfigure}%
  \caption{Second case study: CDDM results on the sequences of four tasks.}
  \label{fig:rve_4tasks}
\end{figure}

However, in an attempt to explore and report on the limitations of the presented method, we also investigated what occurs when we add task D (see Figure \ref{fig:rve_4tasks}) where an RVE has much larger yield stress (see Table \ref{tab:rves_tasks_configurations}, i.e. where the yield stress of the matrix becomes 3.0 MPa instead of 0.5 MPa as in the other tasks. In this case, the plastic response of the RVE is delayed, and we noticed that when Task D is learned first then there would be no advantage in learning cooperatively -- see Figure \ref{fig:rve_4tasks}b. However, if the ordering is different, as shown in Figure \ref{fig:rve_4tasks}a, then the proposed cooperative model is still better. We think it is important to be clear that there might be situations in which the ordering of tasks actually leads to difficulties in learning cooperatively.

\section{Discussion}\label{sec:discussion}

We want to highlight that all the above-mentioned results for both case studies are robust to the hyperparameter choice. For simplicity of presenting the previous results, they refer to a particular architecture configuration. However, in this section we elaborate on the robustness of the CDDM approach presented herein by considering different GRU architectures, varying the number of units and hidden state size. In addition, we elaborate on the knowledge transfer effect and analyze the portion of parameters shared between tasks and the ones dedicated only to one task.

\subsection{Architectures comparison}\label{subsec:arch_comparison}

\begin{figure}[ht!]
	\centering
		\includegraphics[width=\textwidth]{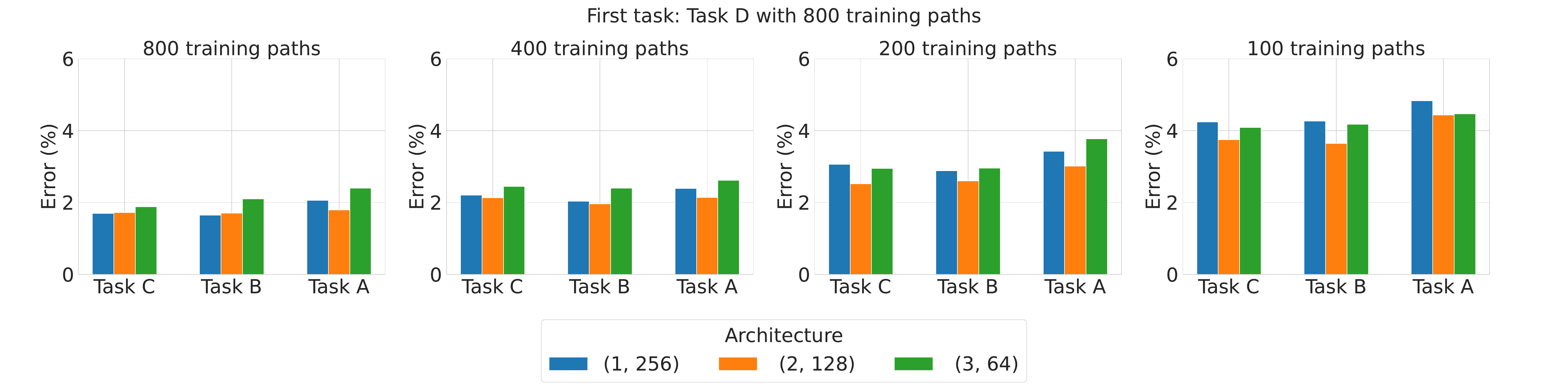}
	  \caption{First case study: architectures comparison on ordering 4.}\label{fig:arch_comparison}
\end{figure}

In this section, we explore how CDDM depends on different GRU architecture hyperparameters such as the number of units or hidden states size. To begin, we consider the first learning case study (Section \ref{sec:plates}) with the corresponding four tasks: the first task is trained with 800 training paths and for the second task we vary the number of training points from 800 to 50. For GRU, we change the number of units from 1 to 3 and consider the hidden state size equal to 64, 128, and 256. Corresponding numbers of learnable parameters are shown in Table \ref{tb:arch_comparison}. In Figure \ref{fig:arch_comparison}, we compare these three network configurations. Overall, we observe similar performance for all of these architectures with insignificant differences. From the figure, we observe that all architectures give us similar results, therefore CDDM is not limited to some special network configuration.

\begin{table}[ht!]
    \centering
    \caption{Architectures comparison.}\label{tb:arch_comparison}
    \begin{tabular}{p{2in} p{1in} p{1in} l}
    \toprule
      Architecture & (1, 256) & (2, 128) & (3, 64) \\
    \midrule
      The number of parameters & 201K & 151K & 63K \\
    \bottomrule
    \end{tabular}
\end{table}

We also want to note that smaller networks (fewer parameters) do not predict the material behavior better when considering the conventional training scenario (non-cooperative). To illustrate this, consider the GRU with 1 cell and the hidden state size of 64 which results in 13K parameters. We train this network on the second case study data and consider average error on tasks A -- D using 800, 400, 200, 100, 50 and 25 training paths. The test errors for these cases are 2.09\%, 2.74\%, 3.59\%, 5.2\%, 7.73\%, 12.9\%. On the other hand, we consider the model with 2 cells and the hidden size of 128 (151K parameters, see Table 4), for which the test errors are 2.01\%, 2.54\%, 3.35\%, 4.84\%, 7.40\%, 11.34\%. As can be seen, the larger model has lower test error.

\subsection{Knowledge transfer}\label{subsec:knowledge_transfer}

\begin{figure}[ht!]
    \begin{subfigure}[b]{0.48\textwidth}
    	\centering
    		\includegraphics[width=\textwidth]{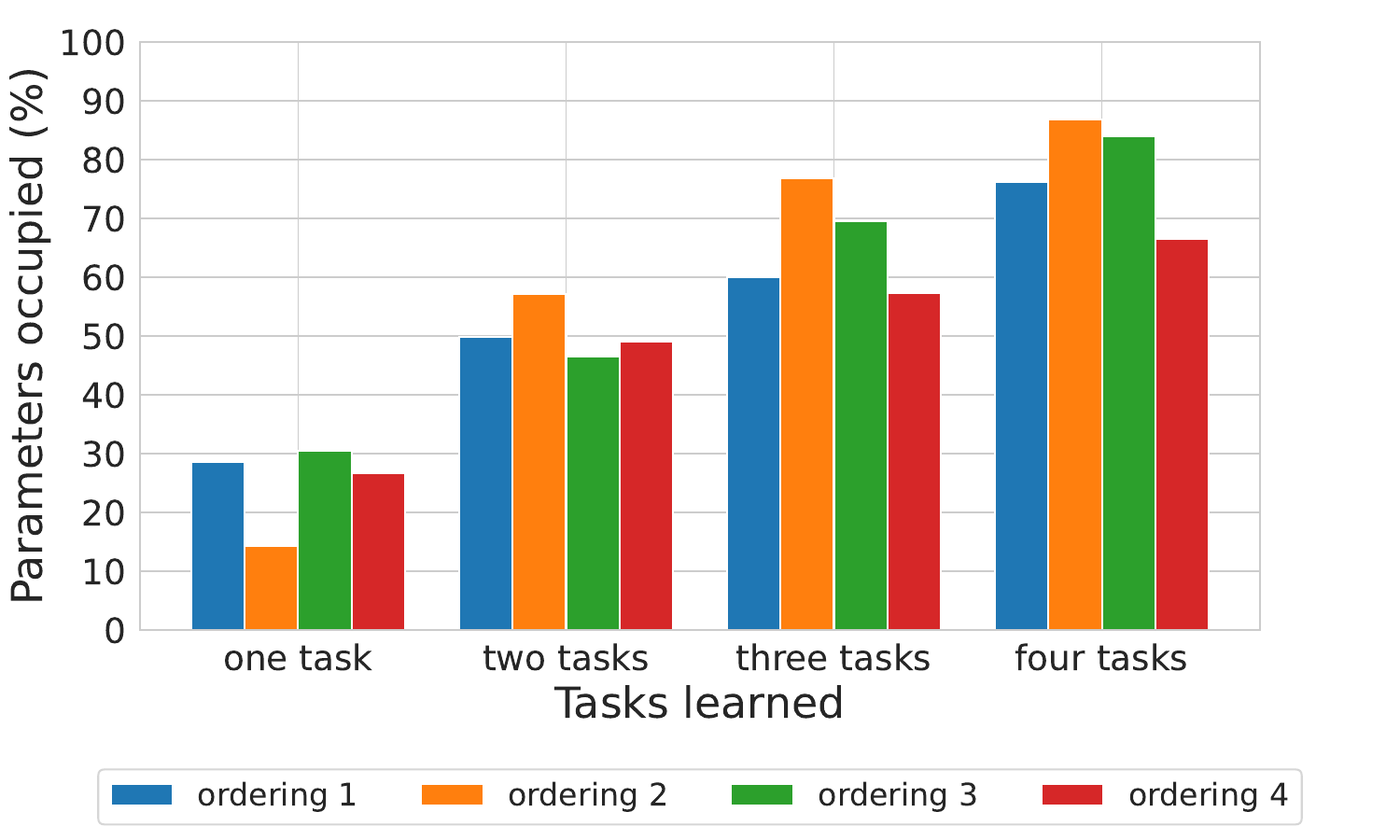}
    	  \caption{Occupied parameters.}
    \end{subfigure}%
    \vspace{1em}
    \begin{subfigure}[b]{0.48\textwidth}
    	\centering
    		\includegraphics[width=\textwidth]{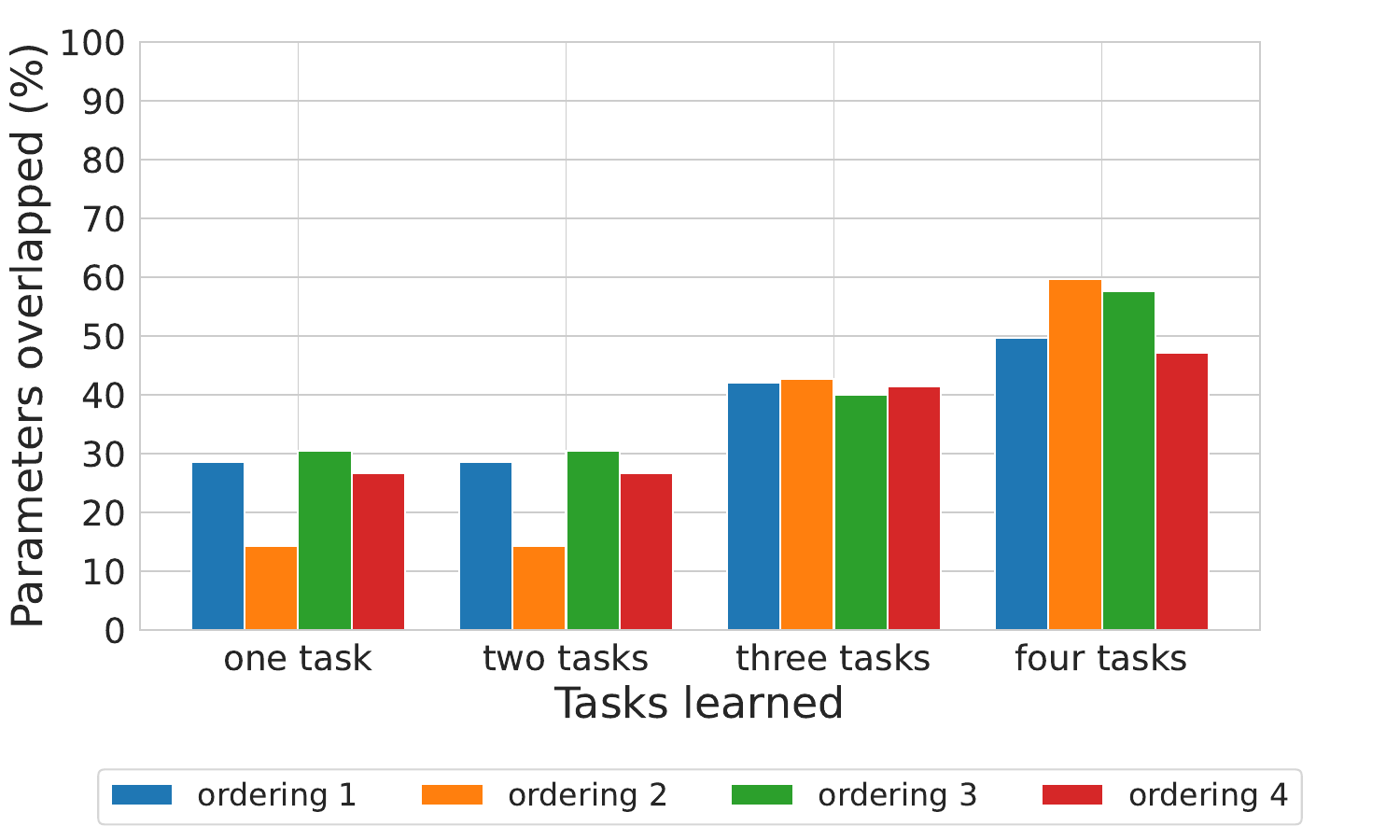}
    	  \caption{Shared parameters.}
    \end{subfigure}
    \caption{First case study: subnetworks analysis: occupied and shared parameters.}
    \label{fig:subnetworks_analysis}
\end{figure}

A crucial characteristic necessary to establishing the CDDM paradigm is the robustness to different task orderings, both in terms of prediction error and the number of parameters used for every task. Focusing on the first case study, we show the relation between subnetworks in the case where the first task is trained with 800 paths and all the other tasks with 200 training paths. In Figure \ref{fig:subnetworks_analysis}a, we show which percentage of the total number of parameters is occupied after a new task is learned. We compare these percentages across four orderings, and in general, we observe the consistency in the number of used parameters with insignificant differences. This means that the ordering of tasks has a negligible effect on how many parameters are occupied in the end.

At the same time, in Figure \ref{fig:subnetworks_analysis}b, we demonstrate the percentage of shared connections while the model learns a new task. So, for instance, we observe in ordering 3 that up to 60\% of parameters are assigned to more than one subnetwork when all the tasks are learned. Overall, at least 45\% of the parameters are shared between multiple tasks without a negative impact on model performance. From the figure, it is clear that the changes in the numbers of shared parameters are consistent across all orderings, illustrating robustness to different task sequences. Nevertheless, after experimenting with many tasks and considering two different case studies, we were able to find one task ordering for Case Study 2 where the cooperative data-driven modeling process was not beneficial when compared to learning all tasks separately (recall Figure \ref{fig:rve_data_comparison}b).

\section{Conclusion and Future Directions}\label{sec:conclusion}

This work introduces the concept of continual learning and the notion of cooperative data-driven modeling. We focus on solid mechanics applications by considering two case studies involving history-dependent plasticity problems. To the best of our knowledge, this is the first example of the application of continual learning in Mechanics and among the first in Engineering applications. We demonstrate that a recurrent neural network can sequentially learn multiple tasks, without replaying data from previous tasks and without forgetting -- an important distinction when comparing to transfer learning methods, and a key enabler of cooperative modeling. The proposed method is based on creating task-related subnetworks that transfer knowledge from each other by sharing neural connections. This is demonstrated to decrease the number of training data required to learn a new task (in this case, a new material law). The approach is robust to different task orders.

As a final note, the authors share their belief that the proposed cooperative data-driven modeling concept has a lot of potential for future developments. More efficient ways of sharing knowledge or selecting subnetworks (when needed), delaying the premature saturation of the network, and accelerating training are only a few possibilities to improve the proposed strategy. Notwithstanding, the prospect of fostering collaborations across different research communities by taking advantage of a machine learning model from a group and adding new capabilities to it such that it solves a new task using less training data and without forgetting how to perform the original task is an exciting new development.

\section*{Acknowledgments}

Miguel A. Bessa acknowledges the support from the NWO Veni award ‘Artificial intelligence towards a sustainable future: ecodesign of recycled polymers and composites’ (with project number 17260 of the research programme Applied and Engineering Sciences) which is financed by the Dutch Research Council (NWO), The Netherlands. Jiaxiang Yi acknowledges the generous support from the China Scholarship Council (CSC).

\appendix

\section{Summary of the pruning algorithm (NNrelief) used by the CP\&S method} \label{appendix:nnrelief}

NNrelief prunes connections with the smallest contribution in the signal of the following neuron. For the incoming signal $\mathbf{X} = \{\mathbf{x}_1, \mathbf{x}_2, \ldots, \mathbf{x}_N\}$ with $N$ data points $\mathbf{x}_n = (x_{n1}, \ldots, x_{nm_{1}}) \in \mathbb{R}^{m_1}$, we compute the importance scores: 
\begin{equation}\label{eq:is}
    s_{ij} (\mathbf{x}_1, \mathbf{x}_2, \ldots, \mathbf{x}_N) = \frac{\overline{\lvert w_{ij} x_{i} \rvert }}{\sum_{k=1}^{m_1} \overline{\lvert w_{kj} x_{k} \rvert}+ \lvert b_j \rvert},
\end{equation}

where $\overline{\lvert w_{ij}x_{i} \rvert} = \frac{1}{N}\sum_{n=1}^N \lvert w_{ij}x_{ni} \rvert $ and $\mathbf{W} = (w_{ij}) \in \mathbb{R}^{m_1 \times m_{2}}$ is a corresponding weight matrix, $\mathbf{b} = (b_1, b_2, \ldots, b_{m_2})^{T} \in \mathbb{R}^{m_{2}}$ is a bias vector. Importance score for the bias of the neuron $j$ is computed by $s_{m_1+1,j} = \frac{\lvert b_j \rvert}{\sum_{k=1}^{m_1} \overline{\lvert w_{kj}x_{k} \rvert} + \lvert b_j \rvert}$. The connections with the smallest importance are pruned. Algorithm \ref{alg:fc_pruning} summarizes the procedure for feedforward layers.

\begin{minipage}{\columnwidth}
    \begin{algorithm}[H]
      \begin{algorithmic}[1]
        \Require network $\mathcal{N}$, training dataset $\mathbf{X}$, pruning hyperparameter $\alpha$.
        \State{$\mathbf{X}^{(0)} \gets \mathbf{X}$}
            \For {every layer $l = 1, \ldots, L$}
                \State{$\mathbf{X}^{(l)}$ $\gets$ layer$\big(\mathbf{X}^{(l-1)}\big)$}
                \For {every neuron $j$ in layer $l$} 
                   \State{Compute importance scores $s_{ij}^{(l)}$ for every incoming connection $w_{ij}$ and bias $b_j$ using Eq. \ref{eq:is}}.
                   \State{$\hat{s}^{(l)}_{ij} \gets \text{Sort}(s^{(l)}_{ij}, order=descending)$}.
                   \State{Find $p_0 = \min\{p : \sum_{i=1}^p \hat{s}^{(l)}_{ij} \ge \alpha \}$}.
                   \State{Prune connections with importance score $s^{(l)}_{ij} < \hat{s}^{(l)}_{p_0j}$}.
                \EndFor
            \EndFor
      \end{algorithmic}
    \caption{Pseudocode for NNrelief}
    \label{alg:fc_pruning}
    \end{algorithm}
\end{minipage}


\bibliographystyle{cas-model2-names}

%\bibliography{References/aleksandr.bib,References/taylan.bib}
\bibliography{references}

\providecommand{\noopsort}[1]{}
\begin{thebibliography}{75}
\expandafter\ifx\csname natexlab\endcsname\relax\def\natexlab#1{#1}\fi
\providecommand{\url}[1]{\texttt{#1}}
\providecommand{\href}[2]{#2}
\providecommand{\path}[1]{#1}
\providecommand{\DOIprefix}{doi:}
\providecommand{\ArXivprefix}{arXiv:}
\providecommand{\URLprefix}{URL: }
\providecommand{\Pubmedprefix}{pmid:}
\providecommand{\doi}[1]{\href{http://dx.doi.org/#1}{\path{#1}}}
\providecommand{\Pubmed}[1]{\href{pmid:#1}{\path{#1}}}
\providecommand{\bibinfo}[2]{#2}
\ifx\xfnm\relax \def\xfnm[#1]{\unskip,\space#1}\fi
%Type = Article
\bibitem[{Abueidda et~al.(2021)Abueidda, Koric, Sobh and Sehitoglu}]{abueidda2021a}
\bibinfo{author}{Abueidda, D.W.}, \bibinfo{author}{Koric, S.}, \bibinfo{author}{Sobh, N.A.}, \bibinfo{author}{Sehitoglu, H.}, \bibinfo{year}{2021}.
\newblock \bibinfo{title}{Deep learning for plasticity and thermo-viscoplasticity}.
\newblock \bibinfo{journal}{International Journal of Plasticity} \bibinfo{volume}{136}, \bibinfo{pages}{102852}.
%Type = Inproceedings
\bibitem[{Aljundi et~al.(2018)Aljundi, Babiloni, Elhoseiny, Rohrbach and Tuytelaars}]{aljundi2018memory}
\bibinfo{author}{Aljundi, R.}, \bibinfo{author}{Babiloni, F.}, \bibinfo{author}{Elhoseiny, M.}, \bibinfo{author}{Rohrbach, M.}, \bibinfo{author}{Tuytelaars, T.}, \bibinfo{year}{2018}.
\newblock \bibinfo{title}{Memory aware synapses: Learning what (not) to forget}, in: \bibinfo{booktitle}{Proceedings of the European Conference on Computer Vision (ECCV)}, pp. \bibinfo{pages}{139--154}.
%Type = Article
\bibitem[{As'~ad et~al.(2022)As'~ad, Avery and Farhat}]{asad2022a}
\bibinfo{author}{As'~ad, F.}, \bibinfo{author}{Avery, P.}, \bibinfo{author}{Farhat, C.}, \bibinfo{year}{2022}.
\newblock \bibinfo{title}{A mechanics-informed artificial neural network approach in data-driven constitutive modeling}.
\newblock \bibinfo{journal}{International Journal for Numerical Methods in Engineering} \bibinfo{volume}{123}, \bibinfo{pages}{2738--2759}.
%Type = Article
\bibitem[{Bengio et~al.(1994)Bengio, Simard and Frasconi}]{bengio1994learning}
\bibinfo{author}{Bengio, Y.}, \bibinfo{author}{Simard, P.}, \bibinfo{author}{Frasconi, P.}, \bibinfo{year}{1994}.
\newblock \bibinfo{title}{Learning long-term dependencies with gradient descent is difficult}.
\newblock \bibinfo{journal}{IEEE transactions on neural networks} \bibinfo{volume}{5}, \bibinfo{pages}{157--166}.
%Type = Article
\bibitem[{Bessa et~al.(2017)Bessa, Bostanabad, Liu, Hu, Apley, Brinson, Chen and Liu}]{bessa2017}
\bibinfo{author}{Bessa, M.}, \bibinfo{author}{Bostanabad, R.}, \bibinfo{author}{Liu, Z.}, \bibinfo{author}{Hu, A.}, \bibinfo{author}{Apley, D.W.}, \bibinfo{author}{Brinson, C.}, \bibinfo{author}{Chen, W.}, \bibinfo{author}{Liu, W.K.}, \bibinfo{year}{2017}.
\newblock \bibinfo{title}{A framework for data-driven analysis of materials under uncertainty: {{Countering}} the curse of dimensionality}.
\newblock \bibinfo{journal}{Computer Methods in Applied Mechanics and Engineering} \bibinfo{volume}{320}, \bibinfo{pages}{633--667}.
\newblock \DOIprefix\doi{10.1016/j.cma.2017.03.037}.
%Type = Article
\bibitem[{Biesialska et~al.(2020)Biesialska, Biesialska and Costa-Jussa}]{biesialska2020continual}
\bibinfo{author}{Biesialska, M.}, \bibinfo{author}{Biesialska, K.}, \bibinfo{author}{Costa-Jussa, M.R.}, \bibinfo{year}{2020}.
\newblock \bibinfo{title}{Continual lifelong learning in natural language processing: A survey}.
\newblock \bibinfo{journal}{arXiv preprint arXiv:2012.09823} .
%Type = Article
\bibitem[{Bonatti et~al.(2022)Bonatti, Berisha and Mohr}]{bonatti2022cp}
\bibinfo{author}{Bonatti, C.}, \bibinfo{author}{Berisha, B.}, \bibinfo{author}{Mohr, D.}, \bibinfo{year}{2022}.
\newblock \bibinfo{title}{From cp-fft to cp-rnn: Recurrent neural network surrogate model of crystal plasticity}.
\newblock \bibinfo{journal}{International Journal of Plasticity} , \bibinfo{pages}{103430}.
%Type = Article
\bibitem[{Capuano and Rimoli(2019)}]{capuano2019smart}
\bibinfo{author}{Capuano, G.}, \bibinfo{author}{Rimoli, J.J.}, \bibinfo{year}{2019}.
\newblock \bibinfo{title}{Smart finite elements: A novel machine learning application}.
\newblock \bibinfo{journal}{Computer Methods in Applied Mechanics and Engineering} \bibinfo{volume}{345}, \bibinfo{pages}{363--381}.
%Type = Article
\bibitem[{Caruana(1994)}]{caruana1994learning}
\bibinfo{author}{Caruana, R.}, \bibinfo{year}{1994}.
\newblock \bibinfo{title}{Learning many related tasks at the same time with backpropagation}.
\newblock \bibinfo{journal}{Advances in neural information processing systems} \bibinfo{volume}{7}.
%Type = Inproceedings
\bibitem[{Chaudhry et~al.(2018)Chaudhry, Dokania, Ajanthan and Torr}]{chaudhry2018riemannian}
\bibinfo{author}{Chaudhry, A.}, \bibinfo{author}{Dokania, P.K.}, \bibinfo{author}{Ajanthan, T.}, \bibinfo{author}{Torr, P.H.}, \bibinfo{year}{2018}.
\newblock \bibinfo{title}{Riemannian walk for incremental learning: Understanding forgetting and intransigence}, in: \bibinfo{booktitle}{Proceedings of the European Conference on Computer Vision (ECCV)}, pp. \bibinfo{pages}{532--547}.
%Type = Article
\bibitem[{Chen(2021)}]{chen2021b}
\bibinfo{author}{Chen, G.}, \bibinfo{year}{2021}.
\newblock \bibinfo{title}{Recurrent neural networks ({{RNNs}}) learn the constitutive law of viscoelasticity}.
\newblock \bibinfo{journal}{Computational Mechanics} , \bibinfo{pages}{11}.
%Type = Article
\bibitem[{Cho et~al.(2014)Cho, Van~Merri{\"e}nboer, Gulcehre, Bahdanau, Bougares, Schwenk and Bengio}]{cho2014learning}
\bibinfo{author}{Cho, K.}, \bibinfo{author}{Van~Merri{\"e}nboer, B.}, \bibinfo{author}{Gulcehre, C.}, \bibinfo{author}{Bahdanau, D.}, \bibinfo{author}{Bougares, F.}, \bibinfo{author}{Schwenk, H.}, \bibinfo{author}{Bengio, Y.}, \bibinfo{year}{2014}.
\newblock \bibinfo{title}{Learning phrase representations using rnn encoder-decoder for statistical machine translation}.
\newblock \bibinfo{journal}{arXiv preprint arXiv:1406.1078} .
%Type = Article
\bibitem[{Dekhovich et~al.(2021)Dekhovich, Tax, Sluiter and Bessa}]{dekhovich2021neural}
\bibinfo{author}{Dekhovich, A.}, \bibinfo{author}{Tax, D.M.}, \bibinfo{author}{Sluiter, M.H.}, \bibinfo{author}{Bessa, M.A.}, \bibinfo{year}{2021}.
\newblock \bibinfo{title}{Neural network relief: a pruning algorithm based on neural activity}.
\newblock \bibinfo{journal}{arXiv preprint arXiv:2109.10795} .
%Type = Article
\bibitem[{Dekhovich et~al.(2022)Dekhovich, Tax, Sluiter and Bessa}]{dekhovich2022continual}
\bibinfo{author}{Dekhovich, A.}, \bibinfo{author}{Tax, D.M.}, \bibinfo{author}{Sluiter, M.H.}, \bibinfo{author}{Bessa, M.A.}, \bibinfo{year}{2022}.
\newblock \bibinfo{title}{Continual prune-and-select: Class-incremental learning with specialized subnetworks}.
\newblock \bibinfo{journal}{arXiv preprint arXiv:2208.04952} .
%Type = Inproceedings
\bibitem[{Douillard et~al.(2020)Douillard, Cord, Ollion, Robert and Valle}]{douillard2020podnet}
\bibinfo{author}{Douillard, A.}, \bibinfo{author}{Cord, M.}, \bibinfo{author}{Ollion, C.}, \bibinfo{author}{Robert, T.}, \bibinfo{author}{Valle, E.}, \bibinfo{year}{2020}.
\newblock \bibinfo{title}{Podnet: Pooled outputs distillation for small-tasks incremental learning}, in: \bibinfo{booktitle}{Computer Vision--ECCV 2020: 16th European Conference, Glasgow, UK, August 23--28, 2020, Proceedings, Part XX 16}, \bibinfo{organization}{Springer}. pp. \bibinfo{pages}{86--102}.
%Type = Article
\bibitem[{Draxl and Scheffler(2018)}]{draxl2018nomad}
\bibinfo{author}{Draxl, C.}, \bibinfo{author}{Scheffler, M.}, \bibinfo{year}{2018}.
\newblock \bibinfo{title}{Nomad: The fair concept for big data-driven materials science}.
\newblock \bibinfo{journal}{Mrs Bulletin} \bibinfo{volume}{43}, \bibinfo{pages}{676--682}.
%Type = Inproceedings
\bibitem[{D{\"u}tting et~al.(2019)D{\"u}tting, Feng, Narasimhan, Parkes and Ravindranath}]{dutting2019optimal}
\bibinfo{author}{D{\"u}tting, P.}, \bibinfo{author}{Feng, Z.}, \bibinfo{author}{Narasimhan, H.}, \bibinfo{author}{Parkes, D.}, \bibinfo{author}{Ravindranath, S.S.}, \bibinfo{year}{2019}.
\newblock \bibinfo{title}{Optimal auctions through deep learning}, in: \bibinfo{booktitle}{International Conference on Machine Learning}, \bibinfo{organization}{PMLR}. pp. \bibinfo{pages}{1706--1715}.
%Type = Article
\bibitem[{French(1999)}]{french1999catastrophic}
\bibinfo{author}{French, R.M.}, \bibinfo{year}{1999}.
\newblock \bibinfo{title}{Catastrophic forgetting in connectionist networks}.
\newblock \bibinfo{journal}{Trends in cognitive sciences} \bibinfo{volume}{3}, \bibinfo{pages}{128--135}.
%Type = Article
\bibitem[{Fuhg and Bouklas(2022)}]{fuhg2022a}
\bibinfo{author}{Fuhg, J.N.}, \bibinfo{author}{Bouklas, N.}, \bibinfo{year}{2022}.
\newblock \bibinfo{title}{The mixed deep energy method for resolving concentration features in finite strain hyperelasticity}.
\newblock \bibinfo{journal}{Journal of Computational Physics} \bibinfo{volume}{451}, \bibinfo{pages}{110839}.
%Type = Article
\bibitem[{Ghaboussi et~al.(1991)Ghaboussi, Garrett and Wu}]{ghaboussi1991a}
\bibinfo{author}{Ghaboussi, J.}, \bibinfo{author}{Garrett, J.H.}, \bibinfo{author}{Wu, X.}, \bibinfo{year}{1991}.
\newblock \bibinfo{title}{Knowledge-based modeling of material behavior with neural networks}.
\newblock \bibinfo{journal}{Journal of Engineering Mechanics} \bibinfo{volume}{117}, \bibinfo{pages}{132--153}.
\newblock \URLprefix \url{https://ascelibrary.org/doi/abs/10.1061/(ASCE)0733-9399(1991)117:1(132)}, \DOIprefix\doi{10.1061/(ASCE)0733-9399(1991)117:1(132)}.
%Type = Article
\bibitem[{Ghavamian and Simone(2019)}]{ghavamian2019a}
\bibinfo{author}{Ghavamian, F.}, \bibinfo{author}{Simone, A.}, \bibinfo{year}{2019}.
\newblock \bibinfo{title}{Accelerating multiscale finite element simulations of history-dependent materials using a recurrent neural network}.
\newblock \bibinfo{journal}{Computer Methods in Applied Mechanics and Engineering} \bibinfo{volume}{357}, \bibinfo{pages}{112594}.
%Type = Article
\bibitem[{Golkar et~al.(2019)Golkar, Kagan and Cho}]{golkar2019continual}
\bibinfo{author}{Golkar, S.}, \bibinfo{author}{Kagan, M.}, \bibinfo{author}{Cho, K.}, \bibinfo{year}{2019}.
\newblock \bibinfo{title}{Continual learning via neural pruning}.
\newblock \bibinfo{journal}{arXiv preprint arXiv:1903.04476} .
%Type = Article
\bibitem[{Goodfellow et~al.(2013)Goodfellow, Mirza, Xiao, Courville and Bengio}]{goodfellow2013empirical}
\bibinfo{author}{Goodfellow, I.J.}, \bibinfo{author}{Mirza, M.}, \bibinfo{author}{Xiao, D.}, \bibinfo{author}{Courville, A.}, \bibinfo{author}{Bengio, Y.}, \bibinfo{year}{2013}.
\newblock \bibinfo{title}{An empirical investigation of catastrophic forgetting in gradient-based neural networks}.
\newblock \bibinfo{journal}{arXiv preprint arXiv:1312.6211} .
%Type = Article
\bibitem[{Han et~al.(2015)Han, Pool, Tran and Dally}]{han2015learning}
\bibinfo{author}{Han, S.}, \bibinfo{author}{Pool, J.}, \bibinfo{author}{Tran, J.}, \bibinfo{author}{Dally, W.J.}, \bibinfo{year}{2015}.
\newblock \bibinfo{title}{Learning both weights and connections for efficient neural networks}.
\newblock \bibinfo{journal}{arXiv preprint arXiv:1506.02626} .
%Type = Article
\bibitem[{Hochreiter(1991)}]{hochreiter1991untersuchungen}
\bibinfo{author}{Hochreiter, S.}, \bibinfo{year}{1991}.
\newblock \bibinfo{title}{Untersuchungen zu dynamischen neuronalen netzen}.
\newblock \bibinfo{journal}{Diploma, Technische Universit{\"a}t M{\"u}nchen} \bibinfo{volume}{91}.
%Type = Article
\bibitem[{Hochreiter and Schmidhuber(1997)}]{hochreiter1997long}
\bibinfo{author}{Hochreiter, S.}, \bibinfo{author}{Schmidhuber, J.}, \bibinfo{year}{1997}.
\newblock \bibinfo{title}{Long short-term memory}.
\newblock \bibinfo{journal}{Neural computation} \bibinfo{volume}{9}, \bibinfo{pages}{1735--1780}.
%Type = Article
\bibitem[{Hu et~al.(2016)Hu, Peng, Tai and Tang}]{hu2016network}
\bibinfo{author}{Hu, H.}, \bibinfo{author}{Peng, R.}, \bibinfo{author}{Tai, Y.W.}, \bibinfo{author}{Tang, C.K.}, \bibinfo{year}{2016}.
\newblock \bibinfo{title}{Network trimming: A data-driven neuron pruning approach towards efficient deep architectures}.
\newblock \bibinfo{journal}{arXiv preprint arXiv:1607.03250} .
%Type = Article
\bibitem[{Ibanez et~al.(2018)Ibanez, Abisset-Chavanne, Aguado, Gonzalez, Cueto and Chinesta}]{ibanez2018a}
\bibinfo{author}{Ibanez, R.}, \bibinfo{author}{Abisset-Chavanne, E.}, \bibinfo{author}{Aguado, J.V.}, \bibinfo{author}{Gonzalez, D.}, \bibinfo{author}{Cueto, E.}, \bibinfo{author}{Chinesta, F.}, \bibinfo{year}{2018}.
\newblock \bibinfo{title}{A manifold learning approach to data-driven computational elasticity and inelasticity}.
\newblock \bibinfo{journal}{Archives of Computational Methods in Engineering} \bibinfo{volume}{25}, \bibinfo{pages}{47--57}.
%Type = Misc
\bibitem[{Jacobsen et~al.(2020)Jacobsen, de~Miranda~Azevedo, Juty, Batista, Coles, Cornet, Courtot, Crosas, Dumontier, Evelo et~al.}]{jacobsen2020fair}
\bibinfo{author}{Jacobsen, A.}, \bibinfo{author}{de~Miranda~Azevedo, R.}, \bibinfo{author}{Juty, N.}, \bibinfo{author}{Batista, D.}, \bibinfo{author}{Coles, S.}, \bibinfo{author}{Cornet, R.}, \bibinfo{author}{Courtot, M.}, \bibinfo{author}{Crosas, M.}, \bibinfo{author}{Dumontier, M.}, \bibinfo{author}{Evelo, C.T.}, et~al., \bibinfo{year}{2020}.
\newblock \bibinfo{title}{Fair principles: interpretations and implementation considerations}.
%Type = Article
\bibitem[{Jones et~al.(2018)Jones, Templeton, Sanders and Ostien}]{jones2018a}
\bibinfo{author}{Jones, R.}, \bibinfo{author}{Templeton, J.}, \bibinfo{author}{Sanders, C.}, \bibinfo{author}{Ostien, J.}, \bibinfo{year}{2018}.
\newblock \bibinfo{title}{Machine learning models of plastic flow based on representation theory.}
\newblock \bibinfo{journal}{CMES-Computer Modeling in Engineering \& Sciences} \bibinfo{volume}{117}.
%Type = Inproceedings
\bibitem[{Jordan(1997)}]{jordan1997serial}
\bibinfo{author}{Jordan, M.I.}, \bibinfo{year}{1997}.
\newblock \bibinfo{title}{Serial order: A parallel distributed processing approach}, in: \bibinfo{booktitle}{Advances in psychology}, \bibinfo{publisher}{Elsevier}. pp. \bibinfo{pages}{471--495}.
%Type = Article
\bibitem[{Karniadakis et~al.(2021)Karniadakis, Kevrekidis, Lu, Perdikaris, Wang and Yang}]{karniadakis2021physics}
\bibinfo{author}{Karniadakis, G.E.}, \bibinfo{author}{Kevrekidis, I.G.}, \bibinfo{author}{Lu, L.}, \bibinfo{author}{Perdikaris, P.}, \bibinfo{author}{Wang, S.}, \bibinfo{author}{Yang, L.}, \bibinfo{year}{2021}.
\newblock \bibinfo{title}{Physics-informed machine learning}.
\newblock \bibinfo{journal}{Nature Reviews Physics} \bibinfo{volume}{3}, \bibinfo{pages}{422--440}.
%Type = Article
\bibitem[{Khalil et~al.(2017)Khalil, Dai, Zhang, Dilkina and Song}]{khalil2017learning}
\bibinfo{author}{Khalil, E.}, \bibinfo{author}{Dai, H.}, \bibinfo{author}{Zhang, Y.}, \bibinfo{author}{Dilkina, B.}, \bibinfo{author}{Song, L.}, \bibinfo{year}{2017}.
\newblock \bibinfo{title}{Learning combinatorial optimization algorithms over graphs}.
\newblock \bibinfo{journal}{Advances in neural information processing systems} \bibinfo{volume}{30}.
%Type = Article
\bibitem[{Kingma and Ba(2014)}]{kingma2014adam}
\bibinfo{author}{Kingma, D.P.}, \bibinfo{author}{Ba, J.}, \bibinfo{year}{2014}.
\newblock \bibinfo{title}{Adam: A method for stochastic optimization}.
\newblock \bibinfo{journal}{arXiv preprint arXiv:1412.6980} .
%Type = Misc
\bibitem[{Lejeune and Zhao(2020)}]{lejeune2020a}
\bibinfo{author}{Lejeune, E.}, \bibinfo{author}{Zhao, B.}, \bibinfo{year}{2020}.
\newblock \bibinfo{title}{Exploring the potential of transfer learning for metamodels of heterogeneous material deformation}.
\newblock \href{http://arxiv.org/abs/2010.16260}{\tt arXiv:2010.16260}.
%Type = Article
\bibitem[{Li and Hoiem(2017)}]{li2017learning}
\bibinfo{author}{Li, Z.}, \bibinfo{author}{Hoiem, D.}, \bibinfo{year}{2017}.
\newblock \bibinfo{title}{Learning without forgetting}.
\newblock \bibinfo{journal}{IEEE transactions on pattern analysis and machine intelligence} \bibinfo{volume}{40}, \bibinfo{pages}{2935--2947}.
%Type = Article
\bibitem[{Liu et~al.(2022)Liu, Kovachki, Li, Azizzadenesheli, Anandkumar, Stuart and Bhattacharya}]{liu2022a}
\bibinfo{author}{Liu, B.}, \bibinfo{author}{Kovachki, N.}, \bibinfo{author}{Li, Z.}, \bibinfo{author}{Azizzadenesheli, K.}, \bibinfo{author}{Anandkumar, A.}, \bibinfo{author}{Stuart, A.M.}, \bibinfo{author}{Bhattacharya, K.}, \bibinfo{year}{2022}.
\newblock \bibinfo{title}{A learning-based multiscale method and its application to inelastic impact problems}.
\newblock \bibinfo{journal}{Journal of the Mechanics and Physics of Solids} \bibinfo{volume}{158}, \bibinfo{pages}{104668}.
%Type = Article
\bibitem[{Liu et~al.(2019a)Liu, Wu and Koishi}]{liu2019deep}
\bibinfo{author}{Liu, Z.}, \bibinfo{author}{Wu, C.}, \bibinfo{author}{Koishi, M.}, \bibinfo{year}{2019}a.
\newblock \bibinfo{title}{A deep material network for multiscale topology learning and accelerated nonlinear modeling of heterogeneous materials}.
\newblock \bibinfo{journal}{Computer Methods in Applied Mechanics and Engineering} \bibinfo{volume}{345}, \bibinfo{pages}{1138--1168}.
%Type = Article
\bibitem[{Liu et~al.(2019b)Liu, Wu and Koishi}]{liu2019c}
\bibinfo{author}{Liu, Z.}, \bibinfo{author}{Wu, C.T.}, \bibinfo{author}{Koishi, M.}, \bibinfo{year}{2019}b.
\newblock \bibinfo{title}{Transfer learning of deep material network for seamless structure\textendash property predictions}.
\newblock \bibinfo{journal}{Computational Mechanics} \bibinfo{volume}{64}, \bibinfo{pages}{451--465}.
\newblock \DOIprefix\doi{10.1007/s00466-019-01704-4}.
%Type = Book
\bibitem[{Logg et~al.(2012)Logg, Mardal and Wells}]{logg2012a}
\bibinfo{editor}{Logg, A.}, \bibinfo{editor}{Mardal, K.A.}, \bibinfo{editor}{Wells, G.} (Eds.), \bibinfo{year}{2012}.
\newblock \bibinfo{title}{Automated {{Solution}} of {{Differential Equations}} by the {{Finite Element Method}}}. volume~\bibinfo{volume}{84} of \textit{\bibinfo{series}{Lecture {{Notes}} in {{Computational Science}} and {{Engineering}}}}.
\newblock \bibinfo{publisher}{{Springer}}, \bibinfo{address}{{Berlin, Heidelberg}}.
\newblock \DOIprefix\doi{10.1007/978-3-642-23099-8}.
%Type = Inproceedings
\bibitem[{Mallya et~al.(2018)Mallya, Davis and Lazebnik}]{mallya2018piggyback}
\bibinfo{author}{Mallya, A.}, \bibinfo{author}{Davis, D.}, \bibinfo{author}{Lazebnik, S.}, \bibinfo{year}{2018}.
\newblock \bibinfo{title}{Piggyback: Adapting a single network to multiple tasks by learning to mask weights}, in: \bibinfo{booktitle}{Proceedings of the European Conference on Computer Vision (ECCV)}, pp. \bibinfo{pages}{67--82}.
%Type = Inproceedings
\bibitem[{Mallya and Lazebnik(2018)}]{mallya2018packnet}
\bibinfo{author}{Mallya, A.}, \bibinfo{author}{Lazebnik, S.}, \bibinfo{year}{2018}.
\newblock \bibinfo{title}{Packnet: Adding multiple tasks to a single network by iterative pruning}, in: \bibinfo{booktitle}{Proceedings of the IEEE conference on Computer Vision and Pattern Recognition}, pp. \bibinfo{pages}{7765--7773}.
%Type = Article
\bibitem[{Masana et~al.(2020a)Masana, Liu, Twardowski, Menta, Bagdanov and van~de Weijer}]{masana2020class}
\bibinfo{author}{Masana, M.}, \bibinfo{author}{Liu, X.}, \bibinfo{author}{Twardowski, B.}, \bibinfo{author}{Menta, M.}, \bibinfo{author}{Bagdanov, A.D.}, \bibinfo{author}{van~de Weijer, J.}, \bibinfo{year}{2020}a.
\newblock \bibinfo{title}{Class-incremental learning: survey and performance evaluation on image classification}.
\newblock \bibinfo{journal}{arXiv preprint arXiv:2010.15277} .
%Type = Article
\bibitem[{Masana et~al.(2020b)Masana, Twardowski and Van~de Weijer}]{masana2020orderings}
\bibinfo{author}{Masana, M.}, \bibinfo{author}{Twardowski, B.}, \bibinfo{author}{Van~de Weijer, J.}, \bibinfo{year}{2020}b.
\newblock \bibinfo{title}{On class orderings for incremental learning}.
\newblock \bibinfo{journal}{arXiv preprint arXiv:2007.02145} .
%Type = Article
\bibitem[{Masi et~al.(2021)Masi, Stefanou, Vannucci and Maffi-Berthier}]{masi2021a}
\bibinfo{author}{Masi, F.}, \bibinfo{author}{Stefanou, I.}, \bibinfo{author}{Vannucci, P.}, \bibinfo{author}{Maffi-Berthier, V.}, \bibinfo{year}{2021}.
\newblock \bibinfo{title}{Thermodynamics-based artificial neural networks for constitutive modeling}.
\newblock \bibinfo{journal}{Journal of the Mechanics and Physics of Solids} \bibinfo{volume}{147}, \bibinfo{pages}{104277}.
%Type = Incollection
\bibitem[{McCloskey and Cohen(1989)}]{mccloskey1989catastrophic}
\bibinfo{author}{McCloskey, M.}, \bibinfo{author}{Cohen, N.J.}, \bibinfo{year}{1989}.
\newblock \bibinfo{title}{Catastrophic interference in connectionist networks: The sequential learning problem}, in: \bibinfo{booktitle}{Psychology of learning and motivation}. \bibinfo{publisher}{Elsevier}. volume~\bibinfo{volume}{24}, pp. \bibinfo{pages}{109--165}.
%Type = Inproceedings
\bibitem[{Mikolov et~al.(2010)Mikolov, Karafi{\'a}t, Burget, Cernock{\`y} and Khudanpur}]{mikolov2010recurrent}
\bibinfo{author}{Mikolov, T.}, \bibinfo{author}{Karafi{\'a}t, M.}, \bibinfo{author}{Burget, L.}, \bibinfo{author}{Cernock{\`y}, J.}, \bibinfo{author}{Khudanpur, S.}, \bibinfo{year}{2010}.
\newblock \bibinfo{title}{Recurrent neural network based language model.}, in: \bibinfo{booktitle}{Interspeech}, \bibinfo{organization}{Makuhari}. pp. \bibinfo{pages}{1045--1048}.
%Type = Article
\bibitem[{Mozaffar et~al.(2019)Mozaffar, Bostanabad, Chen, Ehmann, Cao and Bessa}]{mozaffar2019}
\bibinfo{author}{Mozaffar, M.}, \bibinfo{author}{Bostanabad, R.}, \bibinfo{author}{Chen, W.}, \bibinfo{author}{Ehmann, K.}, \bibinfo{author}{Cao, J.}, \bibinfo{author}{Bessa, M.A.}, \bibinfo{year}{2019}.
\newblock \bibinfo{title}{Deep learning predicts path-dependent plasticity}.
\newblock \bibinfo{journal}{Proceedings of the National Academy of Sciences} \bibinfo{volume}{116}, \bibinfo{pages}{26414--26420}.
\newblock \DOIprefix\doi{10.1073/pnas.1911815116}.
%Type = Article
\bibitem[{Nguyen and Keip(2018)}]{nguyen2018a}
\bibinfo{author}{Nguyen, L.T.K.}, \bibinfo{author}{Keip, M.A.}, \bibinfo{year}{2018}.
\newblock \bibinfo{title}{A data-driven approach to nonlinear elasticity}.
\newblock \bibinfo{journal}{Computers \& Structures} \bibinfo{volume}{194}, \bibinfo{pages}{97--115}.
%Type = Article
\bibitem[{Pan and Yang(2010)}]{pan2010a}
\bibinfo{author}{Pan, S.J.}, \bibinfo{author}{Yang, Q.}, \bibinfo{year}{2010}.
\newblock \bibinfo{title}{A {{Survey}} on {{Transfer Learning}}}.
\newblock \bibinfo{journal}{IEEE Transactions on Knowledge and Data Engineering} \bibinfo{volume}{22}, \bibinfo{pages}{1345--1359}.
\newblock \DOIprefix\doi{10.1109/TKDE.2009.191}.
%Type = Article
\bibitem[{Parisi et~al.(2019)Parisi, Kemker, Part, Kanan and Wermter}]{parisi2019continual}
\bibinfo{author}{Parisi, G.I.}, \bibinfo{author}{Kemker, R.}, \bibinfo{author}{Part, J.L.}, \bibinfo{author}{Kanan, C.}, \bibinfo{author}{Wermter, S.}, \bibinfo{year}{2019}.
\newblock \bibinfo{title}{Continual lifelong learning with neural networks: A review}.
\newblock \bibinfo{journal}{Neural Networks} \bibinfo{volume}{113}, \bibinfo{pages}{54--71}.
%Type = Article
\bibitem[{Paszke et~al.(2019)Paszke, Gross, Massa, Lerer, Bradbury, Chanan, Killeen, Lin, Gimelshein, Antiga et~al.}]{paszke2019pytorch}
\bibinfo{author}{Paszke, A.}, \bibinfo{author}{Gross, S.}, \bibinfo{author}{Massa, F.}, \bibinfo{author}{Lerer, A.}, \bibinfo{author}{Bradbury, J.}, \bibinfo{author}{Chanan, G.}, \bibinfo{author}{Killeen, T.}, \bibinfo{author}{Lin, Z.}, \bibinfo{author}{Gimelshein, N.}, \bibinfo{author}{Antiga, L.}, et~al., \bibinfo{year}{2019}.
\newblock \bibinfo{title}{Pytorch: An imperative style, high-performance deep learning library}.
\newblock \bibinfo{journal}{Advances in neural information processing systems} \bibinfo{volume}{32}.
%Type = Article
\bibitem[{Peng et~al.(2021)Peng, Alber, Buganza~Tepole, Cannon, De, Dura-Bernal, Garikipati, Karniadakis, Lytton, Perdikaris et~al.}]{peng2021multiscale}
\bibinfo{author}{Peng, G.C.}, \bibinfo{author}{Alber, M.}, \bibinfo{author}{Buganza~Tepole, A.}, \bibinfo{author}{Cannon, W.R.}, \bibinfo{author}{De, S.}, \bibinfo{author}{Dura-Bernal, S.}, \bibinfo{author}{Garikipati, K.}, \bibinfo{author}{Karniadakis, G.}, \bibinfo{author}{Lytton, W.W.}, \bibinfo{author}{Perdikaris, P.}, et~al., \bibinfo{year}{2021}.
\newblock \bibinfo{title}{Multiscale modeling meets machine learning: What can we learn?}
\newblock \bibinfo{journal}{Archives of Computational Methods in Engineering} \bibinfo{volume}{28}, \bibinfo{pages}{1017--1037}.
%Type = Inproceedings
\bibitem[{Ramanujan et~al.(2020)Ramanujan, Wortsman, Kembhavi, Farhadi and Rastegari}]{ramanujan2020s}
\bibinfo{author}{Ramanujan, V.}, \bibinfo{author}{Wortsman, M.}, \bibinfo{author}{Kembhavi, A.}, \bibinfo{author}{Farhadi, A.}, \bibinfo{author}{Rastegari, M.}, \bibinfo{year}{2020}.
\newblock \bibinfo{title}{What's hidden in a randomly weighted neural network?}, in: \bibinfo{booktitle}{Proceedings of the IEEE/CVF Conference on Computer Vision and Pattern Recognition}, pp. \bibinfo{pages}{11893--11902}.
%Type = Inproceedings
\bibitem[{Rebuffi et~al.(2017)Rebuffi, Kolesnikov, Sperl and Lampert}]{rebuffi2017icarl}
\bibinfo{author}{Rebuffi, S.A.}, \bibinfo{author}{Kolesnikov, A.}, \bibinfo{author}{Sperl, G.}, \bibinfo{author}{Lampert, C.H.}, \bibinfo{year}{2017}.
\newblock \bibinfo{title}{icarl: Incremental classifier and representation learning}, in: \bibinfo{booktitle}{Proceedings of the IEEE conference on Computer Vision and Pattern Recognition}, pp. \bibinfo{pages}{2001--2010}.
%Type = Article
\bibitem[{Rocha et~al.(2021)Rocha, Kerfriden and {van der Meer}}]{rocha2021a}
\bibinfo{author}{Rocha, I.}, \bibinfo{author}{Kerfriden, P.}, \bibinfo{author}{{van der Meer}, F.}, \bibinfo{year}{2021}.
\newblock \bibinfo{title}{On-the-fly construction of surrogate constitutive models for concurrent multiscale mechanical analysis through probabilistic machine learning}.
\newblock \bibinfo{journal}{Journal of Computational Physics: X} \bibinfo{volume}{9}, \bibinfo{pages}{100083}.
\newblock \URLprefix \url{https://www.sciencedirect.com/science/article/pii/S2590055220300354}, \DOIprefix\doi{https://doi.org/10.1016/j.jcpx.2020.100083}.
%Type = Techreport
\bibitem[{Rumelhart et~al.(1985)Rumelhart, Hinton and Williams}]{rumelhart1985learning}
\bibinfo{author}{Rumelhart, D.E.}, \bibinfo{author}{Hinton, G.E.}, \bibinfo{author}{Williams, R.J.}, \bibinfo{year}{1985}.
\newblock \bibinfo{title}{Learning internal representations by error propagation}.
\newblock \bibinfo{type}{Technical Report}. California Univ San Diego La Jolla Inst for Cognitive Science.
%Type = Article
\bibitem[{Saidi et~al.(2022)Saidi, Pirgazi, Sanjari, Tamimi, Mohammadi, B{\'e}land, Daymond and Tamblyn}]{saidi2022deep}
\bibinfo{author}{Saidi, P.}, \bibinfo{author}{Pirgazi, H.}, \bibinfo{author}{Sanjari, M.}, \bibinfo{author}{Tamimi, S.}, \bibinfo{author}{Mohammadi, M.}, \bibinfo{author}{B{\'e}land, L.K.}, \bibinfo{author}{Daymond, M.R.}, \bibinfo{author}{Tamblyn, I.}, \bibinfo{year}{2022}.
\newblock \bibinfo{title}{Deep learning and crystal plasticity: A preconditioning approach for accurate orientation evolution prediction}.
\newblock \bibinfo{journal}{Computer Methods in Applied Mechanics and Engineering} \bibinfo{volume}{389}, \bibinfo{pages}{114392}.
%Type = Inproceedings
\bibitem[{Sak et~al.(2014)Sak, Senior and Beaufays}]{sak2014long}
\bibinfo{author}{Sak, H.}, \bibinfo{author}{Senior, A.W.}, \bibinfo{author}{Beaufays, F.}, \bibinfo{year}{2014}.
\newblock \bibinfo{title}{Long short-term memory recurrent neural network architectures for large scale acoustic modeling}, in: \bibinfo{booktitle}{Fifteenth Annual Conference of the International Speech Communication Association}.
%Type = Article
\bibitem[{Schmidt et~al.(2019)Schmidt, Marques, Botti and Marques}]{schmidt2019recent}
\bibinfo{author}{Schmidt, J.}, \bibinfo{author}{Marques, M.R.}, \bibinfo{author}{Botti, S.}, \bibinfo{author}{Marques, M.A.}, \bibinfo{year}{2019}.
\newblock \bibinfo{title}{Recent advances and applications of machine learning in solid-state materials science}.
\newblock \bibinfo{journal}{npj Computational Materials} \bibinfo{volume}{5}, \bibinfo{pages}{1--36}.
%Type = Incollection
\bibitem[{Shanmuganathan(2016)}]{shanmuganathan2016artificial}
\bibinfo{author}{Shanmuganathan, S.}, \bibinfo{year}{2016}.
\newblock \bibinfo{title}{Artificial neural network modelling: An introduction}, in: \bibinfo{booktitle}{Artificial neural network modelling}. \bibinfo{publisher}{Springer}, pp. \bibinfo{pages}{1--14}.
%Type = Book
\bibitem[{Smith(2009)}]{abaqus}
\bibinfo{author}{Smith, M.}, \bibinfo{year}{2009}.
\newblock \bibinfo{title}{ABAQUS/Standard User's Manual, Version 6.9}.
\newblock \bibinfo{publisher}{Dassault Syst{\`e}mes Simulia Corp}, \bibinfo{address}{United States}.
%Type = Article
\bibitem[{Sokar et~al.(2021)Sokar, Mocanu and Pechenizkiy}]{sokar2021spacenet}
\bibinfo{author}{Sokar, G.}, \bibinfo{author}{Mocanu, D.C.}, \bibinfo{author}{Pechenizkiy, M.}, \bibinfo{year}{2021}.
\newblock \bibinfo{title}{Spacenet: Make free space for continual learning}.
\newblock \bibinfo{journal}{Neurocomputing} \bibinfo{volume}{439}, \bibinfo{pages}{1--11}.
%Type = Inproceedings
\bibitem[{Tan et~al.(2018)Tan, Sun, Kong, Zhang, Yang and Liu}]{tan2018survey}
\bibinfo{author}{Tan, C.}, \bibinfo{author}{Sun, F.}, \bibinfo{author}{Kong, T.}, \bibinfo{author}{Zhang, W.}, \bibinfo{author}{Yang, C.}, \bibinfo{author}{Liu, C.}, \bibinfo{year}{2018}.
\newblock \bibinfo{title}{A survey on deep transfer learning}, in: \bibinfo{booktitle}{International conference on artificial neural networks}, \bibinfo{organization}{Springer}. pp. \bibinfo{pages}{270--279}.
%Type = Article
\bibitem[{Thakolkaran et~al.(2022)Thakolkaran, Joshi, Zheng, Flaschel, De~Lorenzis and Kumar}]{thakolkaran2022nn}
\bibinfo{author}{Thakolkaran, P.}, \bibinfo{author}{Joshi, A.}, \bibinfo{author}{Zheng, Y.}, \bibinfo{author}{Flaschel, M.}, \bibinfo{author}{De~Lorenzis, L.}, \bibinfo{author}{Kumar, S.}, \bibinfo{year}{2022}.
\newblock \bibinfo{title}{Nn-euclid: deep-learning hyperelasticity without stress data}.
\newblock \bibinfo{journal}{arXiv preprint arXiv:2205.06664} .
%Type = Book
\bibitem[{Thrun and Pratt(1998)}]{thrun1998}
\bibinfo{editor}{Thrun, S.}, \bibinfo{editor}{Pratt, L.} (Eds.), \bibinfo{year}{1998}.
\newblock \bibinfo{title}{Learning to {{Learn}}}.
\newblock \bibinfo{publisher}{{Springer US}}, \bibinfo{address}{{Boston, MA}}.
\newblock \DOIprefix\doi{10.1007/978-1-4615-5529-2}.
%Type = Article
\bibitem[{Vilalta and Drissi(2002)}]{vilalta2002perspective}
\bibinfo{author}{Vilalta, R.}, \bibinfo{author}{Drissi, Y.}, \bibinfo{year}{2002}.
\newblock \bibinfo{title}{A perspective view and survey of meta-learning}.
\newblock \bibinfo{journal}{Artificial intelligence review} \bibinfo{volume}{18}, \bibinfo{pages}{77--95}.
%Type = Article
\bibitem[{Vlassis et~al.(2020)Vlassis, Ma and Sun}]{vlassis2020a}
\bibinfo{author}{Vlassis, N.N.}, \bibinfo{author}{Ma, R.}, \bibinfo{author}{Sun, W.}, \bibinfo{year}{2020}.
\newblock \bibinfo{title}{Geometric deep learning for computational mechanics part i: Anisotropic hyperelasticity}.
\newblock \bibinfo{journal}{Computer Methods in Applied Mechanics and Engineering} \bibinfo{volume}{371}, \bibinfo{pages}{113299}.
%Type = Article
\bibitem[{Wilkinson et~al.(2016)Wilkinson, Dumontier, Aalbersberg, Appleton, Axton, Baak, Blomberg, Boiten, da~Silva~Santos, Bourne et~al.}]{wilkinson2016fair}
\bibinfo{author}{Wilkinson, M.D.}, \bibinfo{author}{Dumontier, M.}, \bibinfo{author}{Aalbersberg, I.J.}, \bibinfo{author}{Appleton, G.}, \bibinfo{author}{Axton, M.}, \bibinfo{author}{Baak, A.}, \bibinfo{author}{Blomberg, N.}, \bibinfo{author}{Boiten, J.W.}, \bibinfo{author}{da~Silva~Santos, L.B.}, \bibinfo{author}{Bourne, P.E.}, et~al., \bibinfo{year}{2016}.
\newblock \bibinfo{title}{The fair guiding principles for scientific data management and stewardship}.
\newblock \bibinfo{journal}{Scientific data} \bibinfo{volume}{3}, \bibinfo{pages}{1--9}.
%Type = Article
\bibitem[{Wortsman et~al.(2020)Wortsman, Ramanujan, Liu, Kembhavi, Rastegari, Yosinski and Farhadi}]{wortsman2020supermasks}
\bibinfo{author}{Wortsman, M.}, \bibinfo{author}{Ramanujan, V.}, \bibinfo{author}{Liu, R.}, \bibinfo{author}{Kembhavi, A.}, \bibinfo{author}{Rastegari, M.}, \bibinfo{author}{Yosinski, J.}, \bibinfo{author}{Farhadi, A.}, \bibinfo{year}{2020}.
\newblock \bibinfo{title}{Supermasks in superposition}.
\newblock \bibinfo{journal}{Advances in Neural Information Processing Systems} .
%Type = Article
\bibitem[{Wu et~al.(2020)Wu, Kilingar, Noels et~al.}]{wu2020recurrent}
\bibinfo{author}{Wu, L.}, \bibinfo{author}{Kilingar, N.G.}, \bibinfo{author}{Noels, L.}, et~al., \bibinfo{year}{2020}.
\newblock \bibinfo{title}{A recurrent neural network-accelerated multi-scale model for elasto-plastic heterogeneous materials subjected to random cyclic and non-proportional loading paths}.
\newblock \bibinfo{journal}{Computer Methods in Applied Mechanics and Engineering} \bibinfo{volume}{369}, \bibinfo{pages}{113234}.
%Type = Inproceedings
\bibitem[{Wu et~al.(2019)Wu, Chen, Wang, Ye, Liu, Guo and Fu}]{wu2019large}
\bibinfo{author}{Wu, Y.}, \bibinfo{author}{Chen, Y.}, \bibinfo{author}{Wang, L.}, \bibinfo{author}{Ye, Y.}, \bibinfo{author}{Liu, Z.}, \bibinfo{author}{Guo, Y.}, \bibinfo{author}{Fu, Y.}, \bibinfo{year}{2019}.
\newblock \bibinfo{title}{Large scale incremental learning}, in: \bibinfo{booktitle}{Proceedings of the IEEE/CVF Conference on Computer Vision and Pattern Recognition}, pp. \bibinfo{pages}{374--382}.
%Type = Article
\bibitem[{Wuest et~al.(2016)Wuest, Weimer, Irgens and Thoben}]{wuest2016machine}
\bibinfo{author}{Wuest, T.}, \bibinfo{author}{Weimer, D.}, \bibinfo{author}{Irgens, C.}, \bibinfo{author}{Thoben, K.D.}, \bibinfo{year}{2016}.
\newblock \bibinfo{title}{Machine learning in manufacturing: advantages, challenges, and applications}.
\newblock \bibinfo{journal}{Production \& Manufacturing Research} \bibinfo{volume}{4}, \bibinfo{pages}{23--45}.
%Type = Inproceedings
\bibitem[{Zenke et~al.(2017)Zenke, Poole and Ganguli}]{zenke2017continual}
\bibinfo{author}{Zenke, F.}, \bibinfo{author}{Poole, B.}, \bibinfo{author}{Ganguli, S.}, \bibinfo{year}{2017}.
\newblock \bibinfo{title}{Continual learning through synaptic intelligence}, in: \bibinfo{booktitle}{International Conference on Machine Learning}, \bibinfo{organization}{PMLR}. pp. \bibinfo{pages}{3987--3995}.
%Type = Article
\bibitem[{Zhang and Mohr(2020)}]{zhang2020using}
\bibinfo{author}{Zhang, A.}, \bibinfo{author}{Mohr, D.}, \bibinfo{year}{2020}.
\newblock \bibinfo{title}{Using neural networks to represent von mises plasticity with isotropic hardening}.
\newblock \bibinfo{journal}{International Journal of Plasticity} \bibinfo{volume}{132}, \bibinfo{pages}{102732}.

\end{thebibliography}

\newpage




\end{document}